\documentclass[pre]{revtex4}

\usepackage{amsmath}
\usepackage{amsfonts}
\usepackage{hyperref}
\usepackage{graphicx}

\begin{document}

\title{Asymptotic Convergence of Weniger's $\delta$-Transformation for a Class of Superfactorially Divergent Stieltjes Series}

\author{Riccardo Borghi}
\affiliation{Dipartimento di Ingegneria Civile, Informatica e delle Tecnologie Aeronautiche, \\
Universit\`{a} ``Roma Tre'', Via Vito Volterra 62, I-00146 Rome, Italy}

\begin{abstract}
The resummation of superfactorially divergent series represents a significant computational challenge in mathematical physics. 
In the present paper the resummation of a specific class of Stieltjes series characterized by a moment sequence growing as $(2n)!$ 
will be addressed. 
Despite the fact that Carleman's condition is satisfied for these series, the convergence rate of Padé approximants is severely hindered by the logarithmic divergence of the associated Carleman series. Weniger's $\delta$ transformation is proposed as a highly efficient alternative resummation tool. By employing recently established results on the converging factors of superfactorially divergent Stieltjes series, an exact integral representation for the truncation error is obtained. This representation enables the rigorous derivation of the leading-order asymptotic behavior of the transformation error, as well as the estimation of the related convergence rate, for real positive arguments. Numerical experiments strongly support the theoretical findings, suggesting that the $\delta$ transformation offers a robust and computationally efficient framework for decoding this class of wildly divergent expansions
\end{abstract}

\maketitle

\section{Introduction}
\label{Sec:Introduction}

{
Divergent series are a fundamental tool in modern theoretical physics, where perturbation expansions often display a factorial or even superfactorial growth of their coefficients ~\citep{Dyson/1952,Bender/Wu/1969}. In practical terms, such series are useful only when paired with suitable resummation techniques. While Borel summation~\citep{Borel/1899,Borel/1928} and Padé approximants~\citep{Pade/1892} are the standard tools, their efficiency significantly drops for 'violently' divergent expansions~\cite{Killingbeck/1977,Weniger/Cizek/Vinette/1993}. A prime example is given by non-relativistic quantum mechanics, where even superfactorial divergence occurs in specific anharmonic oscillators~\citep{Bender/Wu/1971,Bender/Wu/1973}.}
The specific problem addressed in the present paper concerns a class of Stieltjes series with a moment sequence
$\{\mu_n\}^\infty_{n=0}$ growing like $(2n)!$ instead of $n!$ (e.g., the Euler series).
In principle, such series fall within the reach of Padé approximants, as they satisfy the  Carleman condition, which is sufficient 
for Padé resummability. However, since the divergence of the associated Carleman series is merely logarithmic, the 
convergence of Padé approximants is extremely slow, as shown for instance in Ref.~\cite{Weniger/Cizek/Vinette/1993}. 
Then, the resummation problem is tackled by adopting a particular Levin-type nonlinear sequence
transformation, the so-called  Weniger, or $\delta$ transformation~\cite{Weniger/1989}. 
Levin-type transformations~\cite{Levin/1973,Weniger/1989,Weniger/2004} proved to be particularly effective 
and powerful resummation tools, especially when compared to Pad\'e approximants. 
{
From a computational point of view, the practical implementation of Levin-type nonlinear sequence transformations 
for large-scale physical applications was pioneered by J. Grotendorst \cite{Grotendorst/1989,Grotendorst/1991}, 
who developed specific computer algebra frameworks in Maple and FORTRAN to overcome the numerical 
instabilities inherent to the resummation of violently divergent expansions.
}

Despite the empirical effectiveness of Levin-type transformations, and specifically the $\delta$ transformation, 
rigorous analytical convergence proofs remain notably scarce in the literature. While Padé approximants benefit from a 
well established theoretical framework, especially in the case of Stieltjes series ~\citep{Baker/Graves-Morris/1996}, the 
mathematical foundation of
resummation processes involving Levin-type transformations remains largely underdeveloped.
About ten years ago, a {rigorous asymptotic analysis of} the resummation of the Euler series via Weniger's $\delta$ transformation 
has been proposed~\cite{Borghi/Weniger/2015}. 

In the present paper, starting from the same methodology and strategy developed in Ref.~\cite{Borghi/Weniger/2015},
{ we derive an explicit asymptotic estimate for the transformation error of the $\delta$-resummation of Stieltjes series exhibiting 
$(2n)!$ growth. 
 The technical challenges encountered in the present study are considerably harder than those solved in~\cite{Borghi/Weniger/2015}.
By using some recently published results about the converging factors of the class of Stieltjes series under investigation~\cite{Borghi/2025},  an integral representation for the $\delta$ transformation error is first derived, which provides the analytical basis for the subsequent asymptotic analysis.

The paper is organized as follows: Section~\ref{Sec:StieFunStieSer} provides the necessary mathematical preliminaries on Stieltjes series and Levin-type sequence transformations. In particular, much of the material presented in this section is drawn from my previous works, specifically Ref.~\cite{Borghi/Weniger/2015} and Refs.~\cite{Borghi/2024a,Borghi/2024,Borghi/2025}, in order to provide a sufficient self-consistency to the present paper. Nevertheless, readers are strongly encouraged to consult the original publications cited throughout the text. Section~\ref{Sec:InvFactSeries} introduces the integral representation for the transformation error. The core asymptotic estimate is derived in Sections~\ref{Sec:SuperfactorialDivergentStieltjesSeries} and~\ref{Sub:AsymptoticEstimateDeltaConvergenceRate}, with the most delicate analytical steps detailed within several appendices. 
Finally, Section \ref{Sec:Conclusions} presents numerical evidences aimed at supporting the extension of our asymptotic estimates to complex 
values of Stieltjes series arguments.}


{

\section{Decoding Stieltjes series via $\delta$ transformation}
\label{Sec:StieFunStieSer}

\subsection{Stieltjes Series and Carleman's condition}

Consider a positive measure $\mathrm{d}\mu$ on $[0,\infty]$  with finite moment 
\begin{equation}
\label{Eq:Stieltjes.0.1}
\displaystyle
\mu_m=\int_0^\infty\,t^m\,\mathrm{d}\mu\,,\qquad m \ge 0.
\end{equation}

A {{Stieltjes series}} is the~formal power series
\begin{equation}
\label{Eq:Stieltjes.0.2}
\displaystyle
\sum_{m=0}^\infty\,\dfrac{(-1)^m}{z^{m+1}}\,\mu_m\,,
\end{equation}
which is asymptotic, for~ $z\to \infty$,  
to the {{Stieltjes function}} $f(z)$ 
defined as  
\begin{equation}
\label{Eq:Stieltjes.0.3}
\displaystyle
f(z)\,=\,\int_0^\infty\,\frac{\mathrm{d}\mu}{z+t}\,,\qquad |\arg(z)|<\pi\,,
\end{equation}
which is analytic in the cut complex plane $\mathbb{C}  \backslash (-\infty,0])$~\cite{Henrici/1977,Bender/Orszag/1978}. The most known example of Stieltjes series is the Euler series~\cite{Euler/1755}, characterized by the moment sequence $\{\mu_m=m!\}_{m=0}^\infty$, and asymptotic to the so-called Euler integral, 
\begin{equation}
\label{Eq:Stieltjes.3.1}
\displaystyle
\displaystyle\int_0^\infty\,\frac{\exp(-t)\,\mathrm{d}t}{z+t}\,=\,
\exp(z)\,\Gamma(0,z)\,,\qquad |\arg(z)|<\pi.
\end{equation}

The uniqueness of the decoding process of a Stieltjes series~(\ref{Eq:Stieltjes.0.2})  to retrieve the correct value of the corresponding Stieltjes integral~(\ref{Eq:Stieltjes.0.3}) is guaranteed by the Carleman condition, which is a sufficient condition that requires the following series
\begin{equation}
\label{Eq:Calerman}
\displaystyle
\sum_{m=0}^\infty\,\mu^{-\frac 1{2m}}_m\,,
\end{equation}
to be divergent. However, while the Stieltjes series we are going to investigate satisfy this condition, a logarithmic divergence of the series~(\ref{Eq:Calerman}) implies an extremely slow convergence rate for Padé approximants \cite[Ch. 8]{Bender/Orszag/1978}, as noted in~\cite{Weniger/Cizek/Vinette/1993}. 

\subsection{Decoding Stieltjes series via Levin-type transformations}

Given the Stieltjes series  (\ref{Eq:Stieltjes.0.2}), consider the $n$th-order truncation error
\begin{equation}
  \label{Def_StieFunParSum}
  r_n(z)\, = \, f(z)\,-\,f_n(z)\,,
\end{equation}
where $f_n(z)$ denotes the $n$th-order partial sum,
\begin{equation}
  \label{Def_StieFunParSum.1}
  f_n(z) \, = \, \sum_{m=0}^{n} \, \dfrac{(-1)^{m}}{z^{m+1}} \, \mu_{m}.
\end{equation}
The truncation error can be recast as follows:
\begin{equation}
\label{Eq:Stieltjes.7.0}
\displaystyle
r_n(z)\,=\,
\dfrac{(-1)^{n+1}}{z^{n+1}}\,\mu_{n+1}\,\varphi_{n+1}(z)\,,
\end{equation}
where the quantity 
\begin{equation}
\label{Eq:Stieltjes.7}
\displaystyle
\varphi_{n}(z) \,=\,\frac {1}{\mu_{n}}\int_0^\infty\,t^n\,\frac{\mathrm{d}\mu}{t+z}\,,
\qquad n \in \mathbb{N}_0\, ,  \qquad \vert \arg (z) \vert < \pi \, ,
\end{equation}
will henceforth be called the $n$th-order {\em {converging factor}}~\cite{Airey/1937,Dingle/1973} (Actually,  the definition of the converging factor $\varphi_n$ used here differs from the classical definition by a factor $z$. This has been done for making the subsequent calculations easier). 
Levin-type transformations~\cite{Levin/1973}, including the $\delta$ transformation, are constructed by developing approximation models $\varphi_n^{(k)}(z)$ (with $k > 1$) of $\varphi_n(z)$ such that, 
\begin{equation}
\label{Mod_Seq_Om_ConvFact_ForwDiff.1}
\displaystyle\lim_{k\to\infty}\varphi_n^{(k)}(z)=\varphi_n(z),
\end{equation}
in some limiting sense, thus avoiding the direct numerical evaluation  of the integral in  Eq.~(\ref{Eq:Stieltjes.7}). 

\subsection{$\delta$ transformation via annihilation operators}

The $\delta$ transformation can formally be derived by finding suitable linear operators $\hat{T}_{k}$ able to \emph{annihilate} the converging factor model $\varphi_n^{(k)}(z)$, 
\begin{equation}
  \label{Mod_Seq_Om_ConvFact_2}
\hat{T}_{k} \left\{ \varphi_{n}^{(k)} \right\}  \, = \, 0,
\end{equation}
for fixed $k$ but for all $n \in \mathbb{N}_0$.
For the class of Stieltjes series, powerful annihilation operators can be defined via the $k$th-order finite difference operator together with a specific polynomial weight \cite{Weniger/1989},
\begin{equation}
 \label{Mod_Seq_Om_ConvFact_2.2.1}
\hat{T}_{k}\left\{\cdot\right\} = \Delta^k \left\{ (n+\beta)_{k-1} \cdot\right\},
\end{equation}
where $\beta>0$, $(n+\beta)_{k-1}$ denotes the Pochhammer symbol, and the operator $\Delta^k$ 
is defined by~\citep[Eq.~(25.1.1)]{Olver/Lozier/Boisvert/Clark/2010} 
\begin{equation}
  \label{Delta^k_BinomSum}
  \Delta^{k} g(n) \; = \; (-1)^{k} \, \sum_{j=0}^{k} \, (-1)^j \, 
   {\binom {k} {j}} \, g(n+j)\,,\qquad\qquad k \in \mathbb{N}\,.
\end{equation}
This choice is motivated by the fact that the converging factor $\varphi_n$ of any Stieltjes series 
satisfies the first-order difference equation~\cite{Borghi/2024}
\begin{equation}
\label{Eq:Stieltjes.7.3}
\displaystyle
\varphi_{n+1} \,=\,
\frac{\mu_n}{\mu_{n+1}}\,\left(1\,-\,z\,\varphi_n\right)\,,\qquad n \ge 0\,,
\end{equation}
whose solution can be represented as an inverse factorial series,
\begin{equation}
\label{Eq:FactorialExpansionConvergingFactor.1}
\displaystyle
\varphi_{n}(z) \,=\,\sum_{j=0}^\infty\,\dfrac {a_j}{(n+\beta)_{j}}, \qquad\qquad\qquad n \in \mathbb{N}_0\,,
\end{equation}
with $\{a_k\}^\infty_{k=0}$ being a sequence  \emph{independent} of $n$~\cite{Weniger/2010b,Borghi/Weniger/2015,Borghi/2024}.
By substituting into the remainder definition (\ref{Eq:Stieltjes.7.0}) the $k$th-order approximant $\varphi^{(k)}_{n}(z)$ 
obtained by truncating the inverse series in Eq.~(\ref{Eq:FactorialExpansionConvergingFactor.1}) at the first $(k+1)$ terms,
and on using the annihilation operator in Eq.~(\ref{Mod_Seq_Om_ConvFact_2.2.1}), the Stieltjes integral $f(z)$ can be approximated as follows:
\begin{equation}
  \label{Def_StieFunParSum.1.1.1}
f(z)\,\simeq\,
\delta^{(n)}_k(\beta+1),
\end{equation}
where the sequence $\{\delta^{(n)}_k(\gamma)\}_{k=0}^\infty$, with 
\begin{equation}
  \label{Def_StieFunParSum.1.1.2}
\delta^{(n)}_k(\gamma)\,=\,
\dfrac{\Delta^k\left\{(n+\gamma)_{k-1}\dfrac{f_n(z)}{\Delta f_n(z)}\right\}}{\Delta^k\left\{(n+\gamma)_{k-1}\dfrac{1}{\Delta f_n(z)}\right\}}\,,
\qquad\qquad\qquad \gamma>0, n\in \mathbb{N}_0,
\end{equation}
defines the $\delta$ transformation of the partial sum sequence $\{f_n(z)\}_{n=0}^\infty$~\citep{Olver/Lozier/Boisvert/Clark/2010} [Chapter 3.9(v) Levin's and Weniger's Transformations].

E. J. Weniger first employed his transformation for the evaluation of auxiliary functions in molecular electronic structure
calculations~\citep{Weniger/Steinborn/1989a}. Later, it was successfully used  for the
evaluation of special functions~\citep*{Jentschura/Gies/Valluri/Lamm/Weniger/2002,%
  Jentschura/Loetstedt/2012,Jentschura/Mohr/Soff/Weniger/1999,%
  Weniger/1989,Weniger/1990,Weniger/1992,Weniger/1994a,Weniger/1996d,%
  Weniger/2001,Weniger/2008,Weniger/Cizek/1990}, the summation of
divergent perturbation expansions~\citep*{Cizek/Vinette/Weniger/1991,
  Cizek/Vinette/Weniger/1993a,Cizek/Vinette/Weniger/1993b,%
  Caliceti/Meyer-Hermann/Ribeca/Surzhykov/Jentschura/2007,%
  Cizek/Zamastil/Skala/2003,Jentschura/Becher/Weniger/Soff/2000,%
  Jentschura/Weniger/Soff/2000,Weniger/1990,Weniger/1992,Weniger/1994a,%
  Weniger/1996a,Weniger/1996b,Weniger/1996c,Weniger/1996e,%
  Weniger/1997,Weniger/2004,Weniger/Cizek/Vinette/1991,%
  Weniger/Cizek/Vinette/1993}, as well as for the prediction of unknown perturbation
series coefficients
\citep*{Bender/Weniger/2001,Jentschura/Becher/Weniger/Soff/2000,%
  Jentschura/Weniger/Soff/2000,Weniger/1997}. In the last fifteen years, $\delta$ transformation has also been employed in optics in the study of nonparaxial free-space propagation of optical wavefields~\citep*{Borghi/Santarsiero/2003,Li/Zang/Li/Tian/2009,Li/Zang/Tian/2009,%
  Dai/Li/Zang/Tian/2011,Borghi/Gori/Guattari/Santarsiero/2011}, as well as in the
numerical evaluation of several types of stable and unstable diffraction catastrophes~\citep*{Borghi/2007,
Borghi/2012a}.

}


{
\section{The integral representation for the $\delta$ transformation error}
\label{Sec:InvFactSeries}
%

In the present section, an exact integral representation for the transformation error of the $\delta$ transformation,
defined as $f(z)-\delta^{(n)}_k(\gamma)$, will be derived. Such a representation serves as the analytical foundation for the 
asymptotic analysis developed in the subsequent sections.
On applying $\delta$ transformation to both sides of Eq.~(\ref{Def_StieFunParSum}), we have
\begin{equation}
\label{Eq:WenigerTransformationOnK.2}
\begin{array}{l}
\displaystyle
f(z)\,=\,
\dfrac{\Delta^k\left\{(n+\gamma)_{k-1}\dfrac{f_n}{\Delta f_n}\right\}}{\Delta^k\left\{(n+\gamma)_{k-1}\dfrac{1}{\Delta f_n}\right\}}\,+\,
\dfrac{\Delta^k\left\{(n+\gamma)_{k-1}\dfrac{r_n}{\Delta f_n}\right\}}{\Delta^k\left\{(n+\gamma)_{k-1}\dfrac{1}{\Delta f_n}\right\}}\,,
\end{array}
\end{equation}
from which, on taking Eq.~(\ref{Def_StieFunParSum.1.1.2}) and Eq.~(\ref{Eq:Stieltjes.7.0}) into account, the transformation error takes on the following form:
%
\begin{equation}
\label{Eq:WenigerTransformationOnK.2}
\begin{array}{l}
\displaystyle
f(z)\,-\,\delta^{n}_k(\gamma)\,=\,
z\,\dfrac{\Delta^k\left\{(n+\gamma)_{k-1}\varphi_{n+1}\right\}}{\Delta^k\left\{(n+\gamma)_{k-1}\dfrac{1}{\Delta f_n}\right\}}.
\end{array}
\end{equation}
%
%
The asymptotic analysis carried out in Ref.~\cite{Borghi/Weniger/2015} was based on the following integral representation for the $n$th-order converging factor of the Euler series:
\begin{equation}
\label{Eq:WenigerTransformationOnK.2.1}
\begin{array}{l}
\displaystyle
\varphi_n(z)\,=\,\int_0^1\,t^{n-1}\,\exp\left(z\,-\,\dfrac zt\right)\,\mathrm{d}t\,,
\qquad\qquad\qquad\qquad \mathrm{Re}\{z\}>0\,.
\end{array}
\end{equation}
In~\cite{Borghi/2025}, it has been shown how similar integral representations hold for the converging factor of a whole class of superfactorially divergent Stieltjes series, including those under investigation in the present paper.
Accordingly, in the following we shall assume that 

\begin{equation}
  \label{FS_IntRep}
  \varphi_{n}(z) \,=\, \int_{0}^{1} \, t^{n+\beta-1} \, \Phi (t) \, \mathrm{d} t \, , \qquad n+\beta > 0 \, ,
\end{equation}
where the function $\Phi(t)$ acts as the generating function for the sequence $\{a_k\}_{k=0}^\infty$ in the corresponding inverse factorial expansion~\cite{Borghi/2025}.

%
By substituting from Eq.~(\ref{FS_IntRep}) into the numerator of Eq.~(\ref{Eq:WenigerTransformationOnK.2}), 
{and on setting $\gamma=\beta+1$, }
the numerator can be given the following integral representation (see Appendix~\ref{App:Eq:NumeratoreIntegrale} for details):
\begin{equation}
  \label{Eq:NumeratoreIntegrale}
\begin{array}{l}
\displaystyle
\Delta^k\left\{(n+\gamma)_{k-1}\varphi_{n+1}\right\}\,=\,
(-1)^k\,(n+\gamma)_{k-1} \, 
\int_{0}^{1} \, t^{n+\gamma-1} \,{}_{2} F_{1}
\left(
    \genfrac{}{}{0pt}{}{-k,k+n+\gamma-1} {n+\gamma}; t \right) \, \Phi (t)\, \mathrm{d} t  \, ,
\end{array}
 \end{equation}
where ${}_2F_1(\cdot)$ denotes the Gauss hypergeometric function. Then, on substituting from
Eq.~(\ref{Eq:NumeratoreIntegrale}) into Eq.~(\ref{Eq:WenigerTransformationOnK.2}), the following integral representation for the transformation error is 
obtained:
\begin{equation}
\label{Eq:WenigerTransformationError}
\begin{array}{l}
\displaystyle
f(z)\,-\,\delta^{(n)}_k(\gamma)\,=\,z\,(-1)^k\,(n+\gamma)_{k-1} 
\dfrac{
\displaystyle\int_{0}^{1} \, t^{n+\gamma-1} 
\,{}_{2} F_{1}
\left(
    \genfrac{}{}{0pt}{}{-k,k+n+\gamma-1} {n+\gamma}; t \right) \, \Phi (t)\, \mathrm{d} t}
{\Delta^k\left\{(n+\gamma)_{k-1}\dfrac{1}{\Delta f_n}\right\}},
\end{array}
\end{equation}
which constitutes our principal analytical tool for the subsequent asymptotic analysis.
}

\section{Transformation error for a class of superfactorially divergent Stieltjes series}
\label{Sec:SuperfactorialDivergentStieltjesSeries}

{

{
In the present section, the transformation error representation given  in Eq.~(\ref{Eq:WenigerTransformationError}) 
is applied to a specific class of superfactorially divergent Stieltjes series, characterized by the following moment sequence:
\begin{equation}
\label{Eq:Anharmonic.1}
\mu_m\,=\,\Gamma(2 m\,+\,1\,+\,q)\,,\qquad m \ge 0\,,\qquad  q\in(-1,1).
\end{equation}
The restriction of $q$ is dictated by two reasons: the lower bound ensures the integrability of the measure, 
while the upper bound represents a purely technical condition required to ensure the validity of the asymptotic expansions 
carried out in Appendices~\ref{App:Eq:Thm:AsymptoticsDenominatorDelta.1} and~\ref{App:Eq:Thm:AsymptoticsNumeratorDelta.1}.

Moments in Eq.~(\ref{Eq:Anharmonic.1}) are generated by the measure $\mathrm{d}\mu$
\begin{equation}
\label{Eq:Anharmonic.1.0.4}
\mathrm{d}\mu\,=\,\dfrac{t^{\frac{q-1}2}\,\exp(-\sqrt t)}2\,\mathrm{d}t,
\end{equation}
leading to the  Stieltjes function 
\begin{equation}
\label{Eq:Anharmonic.1.4}
f(z)\,=\,
\int_0^\infty\,\dfrac{t^{\frac{q-1}2}\,\exp(-\sqrt t)}{2(z+t)}\,\mathrm{d}t.
\end{equation}
The closed-form of $f(z)$, utilized later for numerical validation, is derived in Appendix~\ref{App:AnalyticalStieltjes}.

To evaluate the denominator of the expression of the transformation error in Eq.~(\ref{Eq:WenigerTransformationError}), 
we evaluate the explicit form of  $1/\Delta f_n$, 
\begin{equation}
\label{Eq:WenigerTransformationErrorBis.1.1}
\begin{array}{l}
\displaystyle
\dfrac 1{\Delta f_n}\,=\,\dfrac{(-)^{n+1}}{\Gamma(2 n\,+\,3\,+\,q)}\,z^{n+2},
\end{array}
\end{equation}
which, after some algebra, gives
\begin{equation}
\label{Eq:WenigerTransformationErrorBis.3}
\begin{array}{l}
\displaystyle
\Delta^k\left\{(n+\gamma)_{k-1}\,\dfrac{(-)^{n+1}}{\mu_{n+1}}\,z^{n+2}\right\}\,=\,\\
\\
\displaystyle
\,=\,
-z^2\,(-)^n\,(-1)^k\,\dfrac{(n+\gamma)_{k-1}}{\Gamma(2n+q+3)}\,
\,{}_{2} F_{3}\left(
    \genfrac{}{}{0pt}{}{-k,k+n+\gamma-1} {n+\dfrac{q}{2}+\dfrac{3}{2},n+\dfrac{q}{2}+2,
   n+\gamma}; -\dfrac{z}{4} \right)\,,
\end{array}
\end{equation}
and, once substituted into Eq.~(\ref{Eq:WenigerTransformationError}), leads to
\begin{equation}
\label{Eq:WenigerTransformationErrorTer}
\begin{array}{l}
\displaystyle
f(z)\,-\,\delta^{(n)}_k(\gamma)\,=\,\\
\\
\,=\,
\left(-\dfrac 1z\right)^{n+1}\,\Gamma(2n+q+3)\,
\dfrac{
\displaystyle\int_{0}^{1} \, t^{n+\gamma-1} 
\,{}_{2} F_{1}
\left(
    \genfrac{}{}{0pt}{}{-k,k+n+\gamma-1} {n+\gamma}; t \right) \, \Phi (t)\, \mathrm{d} t}
{\displaystyle {}_{2} F_{3}\left(
    \genfrac{}{}{0pt}{}{-k,k+n+\gamma-1} {n+\dfrac{q}{2}+\dfrac{3}{2},n+\dfrac{q}{2}+2,
   n+\gamma}; -\dfrac{z}{4} \right)}.
\end{array}
\end{equation}
As far as the generating function $\Phi(t)$ is concerned, in Ref.~\cite{Borghi/2025} it was proved that,
for the moment sequence~(\ref{Eq:Anharmonic.1}) with $q=0$, 
\begin{equation}
\label{Eq:DifferentialEquation.13.New}
\begin{array}{l}
\displaystyle
\Phi(t)\,=\,\dfrac 1{2\sqrt z}\,\sin\left[\sqrt z\,\left(\dfrac 1{\sqrt{t}}\,-\,1\right)\right]\,,\qquad\qquad\qquad\qquad\qquad\qquad z>0,
\end{array}
\end{equation}
with the parameter $\beta$ in Eq.~(\ref{FS_IntRep}) being set to zero. 

When $q\ne 0$, the function $\Phi(t)$ in Eq.~(\ref{Eq:DifferentialEquation.13.New}) can still be used to represent $\varphi_n(z)$, provided that $\beta=q/2$. To prove this, we substitute from Eq.~(\ref{Eq:Anharmonic.1}) into Eq.~(\ref{Eq:Stieltjes.7.3}), which after some algebra gives
\begin{equation}
\label{Eq:NewRecurrence}
\displaystyle
4\,\left(n+\dfrac q2+1\right)\,\left(n+\dfrac q2+\dfrac12\right)\,\varphi_{n+1}\,+\,z\,\varphi_n \,=\,1,\qquad\qquad n>0,
\end{equation}
and formally coincides with~\cite[Eq.~(52)]{Borghi/2025} after letting $n+q/2=m$. Then, 
on using the procedure described by Eqs. (52)-(58) of~\cite{Borghi/2025}, it is found that, for $z>0$, 
the $n$th-order converging factor of the Stieltjes series in Eq.~(\ref{Eq:Anharmonic.1}) can be expressed as
\begin{equation}
\label{Eq:NewRecurrence.2}
\displaystyle
\varphi_n(z)\,=\,\int_0^1\,t^{n+\frac q2-1}\Phi(t)\,\mathrm{d}t,\qquad\qquad n>0,
\end{equation}
with $\Phi(t)$ given by Eq.~(\ref{Eq:DifferentialEquation.13.New}).

Finally, on adopting the standard choice $n=0$ with $\gamma=1+q/2$, the $k$th-order transformation error 
$\mathcal{E}_k(z,q)\,=\,f(z)\,-\,\delta^{(0)}_k(1+q/2)$ takes the form
\begin{equation}
\label{Eq:WenigerTransformationErrorQuaTer.3.1}
\begin{array}{l}
\displaystyle
\mathcal{E}_k(z,q)\,=\,
\left(-\dfrac 1z\right)\,\Gamma(q+3)\,
\dfrac{
\displaystyle\int_{0}^{1} \, t^{q/2} 
\,{}_{2} F_{1}
\left(
    \genfrac{}{}{0pt}{}{-k,k+q/2} {1+q/2}; t \right) \, \Phi (t)\, \mathrm{d} t}
{\displaystyle {}_{2} F_{3}\left(
    \genfrac{}{}{0pt}{}{-k,k+\dfrac q2} {\dfrac q2+1,\dfrac{q}{2}+\dfrac{3}{2},\dfrac{q}{2}+2}; -\dfrac{z}{4} \right)}\,.
\end{array}
\end{equation}
Equation (\ref{Eq:WenigerTransformationErrorQuaTer.3.1}) represents the starting point for the asymptotic analysis 
presented in the next section.
While a rather similar algebraic representation of the remainder in terms of hypergeometric functions was already derived 
for the Euler series in \cite{Borghi/Weniger/2015} (where a ${}_2F_2$ hypergeometric polynomial was involved), its extension to the present 
class represents one of the main contributions of the present work.
}

\section{Asymptotic estimate of the delta convergence rate}
\label{Sub:AsymptoticEstimateDeltaConvergenceRate}

{

In the present section, the asymptotic behavior, for large values of $k$, of the transformation error defined in 
Eq.~(\ref{Eq:WenigerTransformationErrorQuaTer.3.1}), will be obtained. The analysis is carried out for $z>0$ and $q\in(-1,1)$
and, for the sake of  clarity, the most technical mathematical derivations are detailed in Appendices~\ref{App:Eq:Thm:AsymptoticsDenominatorDelta.1} and~\ref{App:Eq:Thm:AsymptoticsNumeratorDelta.1}.

The $\delta$ transformation is expressed as the ratio of two polynomials,  so that is a meromorphic function
whose poles corresponds to the zeros of the denominator polynomial in Eq.~(\ref{Eq:WenigerTransformationErrorQuaTer.3.1}). 
For the class of Stieltjes series under investigation, it has recently been established~\cite{Slevinsky/2025,Borghi/2026}  that zeros are all real and negative. As a consequence, the $\delta$ transformation correctly simulates the branch cut 
of the Stieltjes function $f(z)$ for all $k\in\mathbb{N}$ and $q\in(-1,1)$.
The transformation error will now be asymptotically estimated within the limit $k\gg 1$. 

The denominator of $\mathcal{E}_k(z,q)$ involves a ${}_2F_3$ hypergeometric polynomial. 
Its asymptotic expansion for $k\to\infty$ is derived in Appendix~\ref{App:Eq:Thm:AsymptoticsDenominatorDelta.1},
where it is found
\begin{equation}
\label{Eq:Thm:AsymptoticsDenominatorDelta.1}
\begin{array}{l}
\displaystyle
{{}_{2} F_{3} \left(
    \genfrac{}{}{0pt}{}{-k, k+\dfrac q2} 
     {\rho_{q}}; -\dfrac z4 \right) \; \sim \;}\\
\\
\sim
\dfrac { 
\Gamma\left(\dfrac q2+1\right)\,\Gamma\left(\dfrac q2+\dfrac 32\right)\,\Gamma\left(\dfrac q2+2\right)}
{2(2\pi)^{3/2}}\,\left(\dfrac z4\right)^{-\frac 38(q+2)}\,k^{-\frac 34(q+2)}\, \exp \left\{2\sqrt 2 k^{1/2} z^{1/4}\right\},
 \quad {k \to \infty}.
\end{array}
\end{equation}
The numerator in Eq.~(\ref{Eq:WenigerTransformationErrorQuaTer.3.1}) requires the asymptotic evaluation of an integral involving the Gauss hypergeometric function ${}_2F_1$. 
As shown in Appendix~\ref{App:Eq:Thm:AsymptoticsNumeratorDelta.1}, for $k\to\infty$, the following approximation holds:
\begin{equation}
\label{Eq:Thm:AsymptoticsNumeratorDelta.1}
\begin{array}{l}
\displaystyle
 {}_{2} F_{1} \left(
   \genfrac{}{}{0pt}{}{-k, k+\dfrac q2}      
    {1+\dfrac q2}; t \right) \; \sim \;\\
\\
\sim\,
\displaystyle
\dfrac {\Gamma\left(\dfrac q2+1\right)}{\sqrt\pi} 
\,k^{-\frac {1+q}2}\,
         \sin^{-\frac{1+q}2} \dfrac\theta 2\,\sqrt{\cos \dfrac\theta 2}
          \cos \left(k \theta - \dfrac{1+q}4\pi\right)\,,\qquad\qquad k \to \infty,
\end{array}
\end{equation}
where $t=\sin^2\dfrac\theta 2$.
}
Then, on substituting from Eqs. (\ref{Eq:DifferentialEquation.13.New}) and (\ref{Eq:Thm:AsymptoticsNumeratorDelta.1}) into the numerator of Eq.~(\ref{Eq:WenigerTransformationErrorQuaTer.3.1}), after changing the integration variable 
$t \in [0,1]$ into the new variable $\theta \in [0,\pi]$, we have
\begin{equation}
  \label{Eq:ExactNumeratorN0Beta1}
\begin{array}{l}
\displaystyle
\displaystyle\int_{0}^{1} \, t^{q/2} 
\,{}_{2} F_{1}
\left(
    \genfrac{}{}{0pt}{}{-k,k+q/2} {1+q/2}; t \right) \, \Phi (t)\, \mathrm{d} t\,\sim\,
\dfrac {\Gamma\left(\dfrac q2+1\right)}{2\,\sqrt z\,\sqrt\pi} 
\,k^{-\frac {1+q}2}\,\\
\\
\times
\displaystyle\int_0^\pi\,
\sin\left[\sqrt z\,\left(\dfrac 1{\sin\dfrac\theta 2}\,-\,1\right)\right]\,
\cos\left(k\theta-\dfrac{1+q}4\pi\right)\,
\sin^{\frac{1+q}2}\dfrac\theta 2\,\cos^{3/2}\dfrac\theta 2\,\mathrm{d}\theta\,,\qquad k\to\infty.
\end{array}
 \end{equation}
An asymptotic estimate of the $\theta$-integral in Eq.~(\ref{Eq:ExactNumeratorN0Beta1}) has been detailed in  Appendix~\ref{App:ExactNumeratorAsymtpt.1},
where it is proved that
\begin{equation}
  \label{Eq:ExactNumeratorAsymtpt.1}
\begin{array}{l}
\displaystyle
\displaystyle\int_0^\pi\,
\sin\left[\sqrt z\,\left(\dfrac 1{\sin\dfrac\theta 2}\,-\,1\right)\right]\,
\cos\left(k\theta-\dfrac{1+q}4\pi\right)\,
\sin^{\frac{1+q}2}\dfrac\theta 2\,\cos^{3/2}\dfrac\theta 2\,\mathrm{d}\theta\,\\
\\
\displaystyle
\sim
\dfrac{\sqrt\pi\, z^{\frac{2+q}8}}{(2k)^{1+\frac q4}}\,\sin\left(2\sqrt{2} z^{1/4} k^{1/2}-\sqrt{z}-\dfrac{q\pi}4\right)
,\qquad\qquad k\to\infty,
\end{array}
 \end{equation}
which, once inserted into Eq.~(\ref{Eq:ExactNumeratorN0Beta1}) gives at once
\begin{equation}
  \label{Eq:ExactNumeratorAsymtpt.1.1}
\begin{array}{l}
\displaystyle
\displaystyle\int_{0}^{1} \, t^{q/2} 
\,{}_{2} F_{1}
\left(
    \genfrac{}{}{0pt}{}{-k,k+q/2} {1+q/2}; t \right) \, \Phi (t)\, \mathrm{d} t\,\sim\,
\,\\
\\
\dfrac{\Gamma\left(\dfrac q2+1\right)}{2^{2+\frac q4}}\, 
\,k^{-\frac 34 (2+q)}\,z^{\frac 18(q-2)}\,\sin\left(2\sqrt{2} z^{1/4} k^{1/2}-\sqrt{z}-\dfrac{q\pi}4\right)\,,\qquad k\to\infty.
\end{array}
 \end{equation}
Finally, on substituting from Eqs. (\ref{Eq:Thm:AsymptoticsDenominatorDelta.1}) and (\ref{Eq:ExactNumeratorAsymtpt.1.1}) into Eq.~(\ref{Eq:WenigerTransformationErrorQuaTer.3.1}), the following asymptotics of the transformation error is eventually obtained:
\begin{equation}
\label{Eq:WenigerTransformationErrorAsymptotics}
\begin{array}{l}
\displaystyle
\mathcal{E}_k(z,q)\,\sim\,
-2\pi\,z^{\frac{q-1}2}\,\exp \left(-2\sqrt 2\, z^{1/4} \, k^{1/2}\right)\,\sin\left(2\sqrt{2}\, z^{1/4}\, k^{1/2}-\sqrt{z}-\dfrac{q\pi}4\right),\quad k\to\infty.
\end{array}
\end{equation}
Equation (\ref{Eq:WenigerTransformationErrorAsymptotics}) constitutes the main result of the present paper.
It predicts an exponential convergence, similarly as it was found in Ref.~\cite{Borghi/Weniger/2015} for the Euler series,
where the transformation error was proportional to $\exp\left(-\frac 92\,z^{1/3}\,k^{2/3}\right)$. 
In the present case, the exponential factor $\exp(- 2\sqrt 2\,z^{1/2}\,k^{1/2})$  reflects an increased difficulty in resumming 
the series, due to its superfactorial growth.

{
While the leading-order estimate in Eq.~(\ref{Eq:WenigerTransformationErrorAsymptotics}) has been derived only for
$z>0$, its validity across the cut complex plane $\mathbb{C}\backslash (-\infty,0]$ is here proposed in the form of an analytical conjecture. In particular, following the approach of Ref. \cite{Borghi/Weniger/2015}, we adopt the principal branch of the fractional powers of $z^{1/4}$. Accordingly, from a formal standpoint, a similar error decay
should persist as long as  the condition  $\mathrm{Re}\{z^{1/4}\}>0$ is met, condition that is satisfied throughout the principal sheet of the associated Riemann surface. 
It must be stressed that the above consideration represents nothing but a conjecture,  whose robustness will be supported only by
numerical experiments.
This conjecture would suggest that the $\delta$ transformation effectively "mimics" the analyticity of the Stieltjes function 
by clustering its poles along the negative real axis, thus providing a stable resummation even in the half-plane 
$\mathrm{Re}\{z\}<0)$, provided the branch cut is avoided. The numerical simulations that will be illustrated in the next section seem to confirm such a conjecture.
}

\section{Numerical Results}
\label{Sec:NumericalResults}

{
The present section provides a numerical validation for the asymptotic convergence rate 
in Eq.~(\ref{Eq:WenigerTransformationErrorAsymptotics}). To this end, the $\delta$ transformation has been implemented 
using high-precision arithmetic (via {\em Wolfram Mathematica} 14.3). In particular,  to ensure stability at high transformation orders 
$(k\le 100)$ and to prevent catastrophic cancellations of digits typical of these nonlinear mappings, the 
$\delta$ transformation has been implemented using arbitrary-precision arithmetic (via Wolfram Mathematica 14.3) with a working precision set to 150 digits.
}

Figure~\ref{Fig:NumericalResults.1} compares the experimental data obtained by using 
the $\delta$ transformation and the asymptotic estimate provided by Eq.~(\ref{Eq:WenigerTransformationErrorAsymptotics}). We refer to the resummation of the Stieltjes series, for the pair $(z,q)=(1,0)$. Open circles represent the modulus of transformation error values $\mathcal{E}_k(1,0)$ obtained through the exact value of $f(z)$ provided by Eq.~(\ref{Eq:Anharmonic.1.4}), explicitly given in Appendix~\ref{App:AnalyticalStieltjes}, and the $\delta$ sequence values $\delta_k^{(0)}(1)$,
with $k \in \{2, 3, 4, \ldots, K\}$, with $K=100$.
The solid curve is the corresponding asymptotic estimate of the transformation error given in Eq.~(\ref{Eq:WenigerTransformationErrorAsymptotics}).
\begin{figure}[ht!]
\includegraphics[width=10cm]{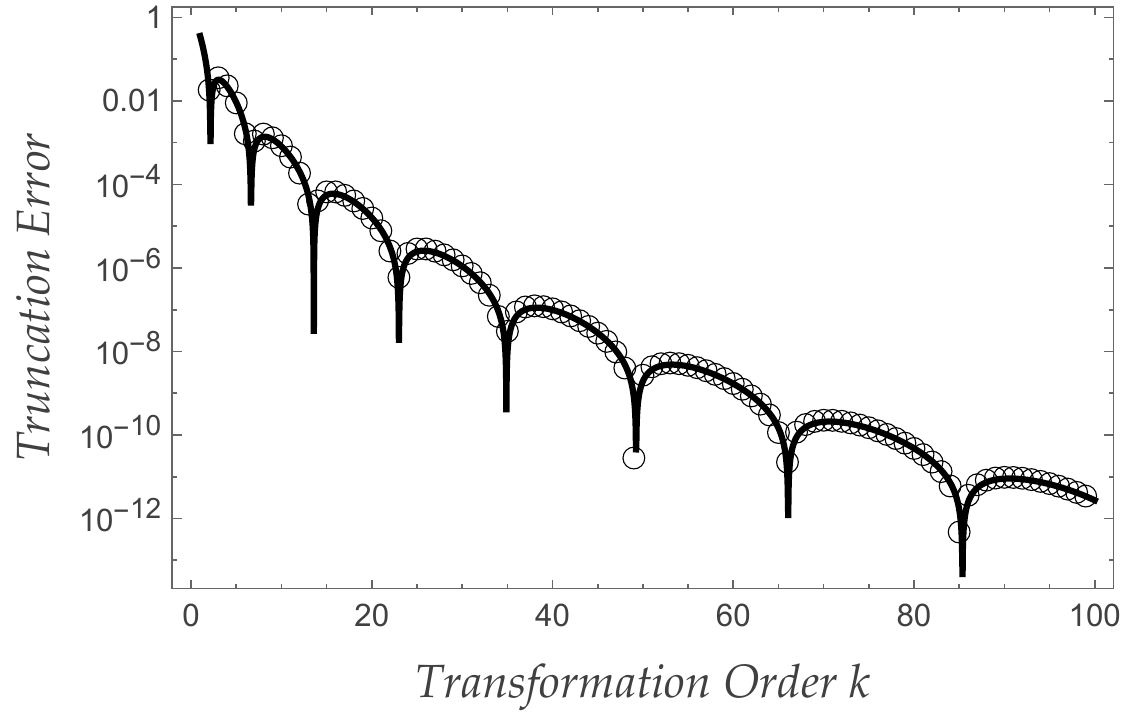}
\caption{Comparison of the observed transformation errors (open circles) vs the transformation order $k$
of the sequence 	$\delta_k^{(0)}(1+q/2)$ for the resummation of the Stieltjes function corresponding 
to $(z,q)=(1,0)$. The solid curve represents the corresponding asymptotic estimate given by Eq.~(\ref{Eq:WenigerTransformationErrorAsymptotics}).
{ The periodic dips observed in the error curves are governed by the trigonometric factor in Eq.~(\ref{Eq:WenigerTransformationErrorAsymptotics}). 
Specifically, for $(z,q)=(1,0)$, the error oscillates as $\sin(2\sqrt 2 \sqrt  k-\pi/4)$, resembling a 'chirp signal' with an instantaneous spatial frequency 
that decreases as $1/\sqrt k$ (such insightful analogy was suggested to me by one of the anonymous Reviewers).} }
\label{Fig:NumericalResults.1}
\end{figure}

The Stieltjes series (\ref{Eq:Stieltjes.0.2}) is asymptotic to the Stieltjes integral 
(\ref{Eq:Stieltjes.0.3}) as $z\to\infty$. Therefore, the argument $z=1$ represents a 
highly challenging summation problem.  The results depicted in  Fig.~\ref{Fig:NumericalResults.1} also 
demonstrate how the asymptotic estimate (\ref{Eq:WenigerTransformationErrorAsymptotics}) 
accurately describes the error behavior even for relatively low values of $k$, well before the 
formal asymptotic regime is reached.

{
In order to numerically validate  the analytical conjecture presented at the end of  Section~\ref{Sub:AsymptoticEstimateDeltaConvergenceRate}, the $\delta$ transformation has been applied also
to complex arguments of the form $z=\exp(\mathrm{i}\phi)$, for several phases $\phi$. 
As shown in Fig.~\ref{Fig:NumericalResults.2}, the agreement remains robust even as the argument $z$ approaches the branch cut ($\phi=9\pi/10$), although a predictable increase in the error magnitude clearly appears.
}
%
%
\begin{figure}[ht!]
  \centering
    \begin{minipage}{6cm}
        \centering
        \includegraphics[width=6cm]{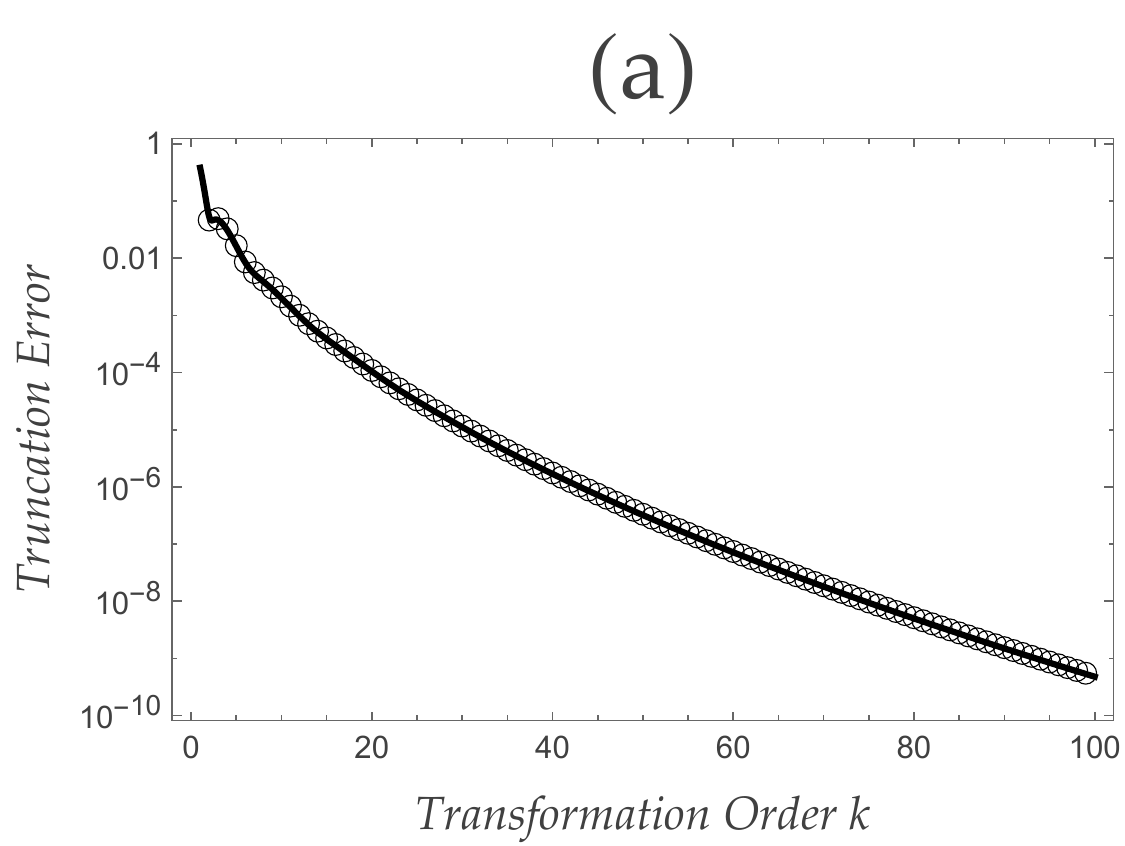}
     \end{minipage}
    \begin{minipage}{6cm}
        \centering
        \includegraphics[width=6cm]{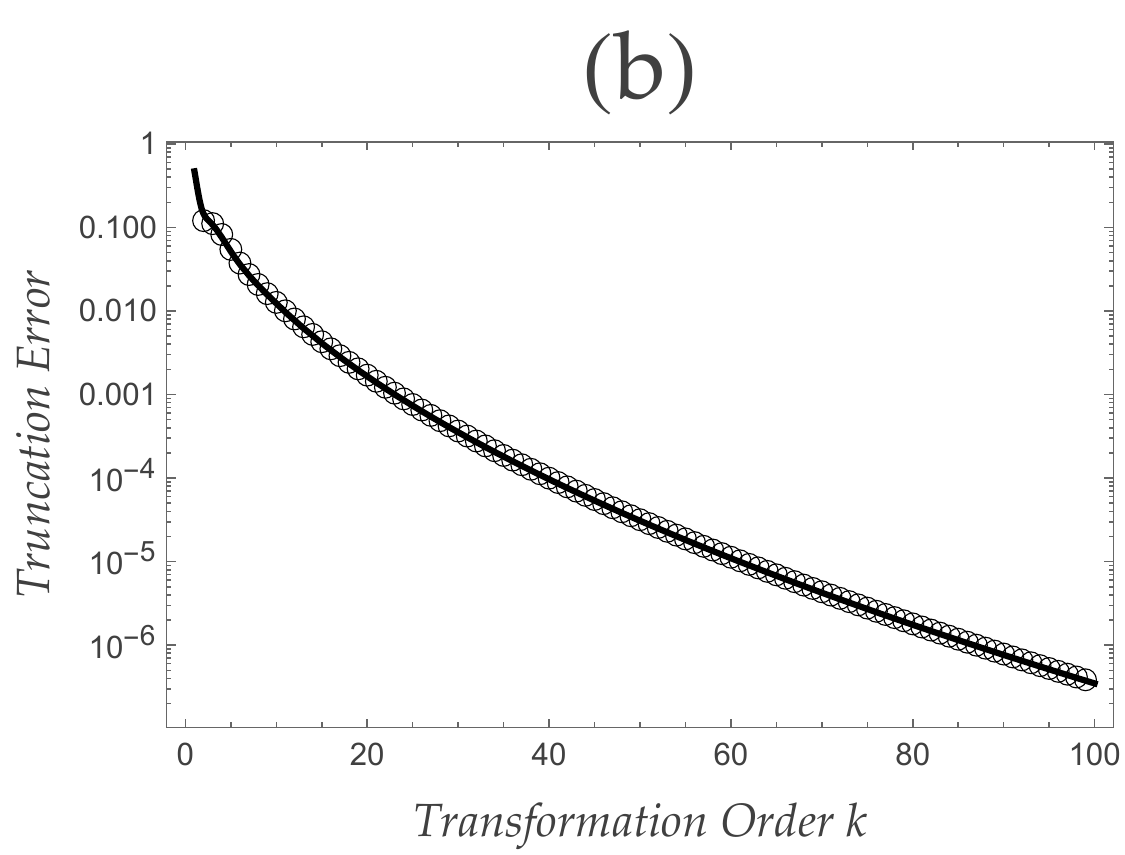}
    \end{minipage}

    \vspace{0.5cm} 

    \begin{minipage}{6cm}
        \centering
        \includegraphics[width=6cm]{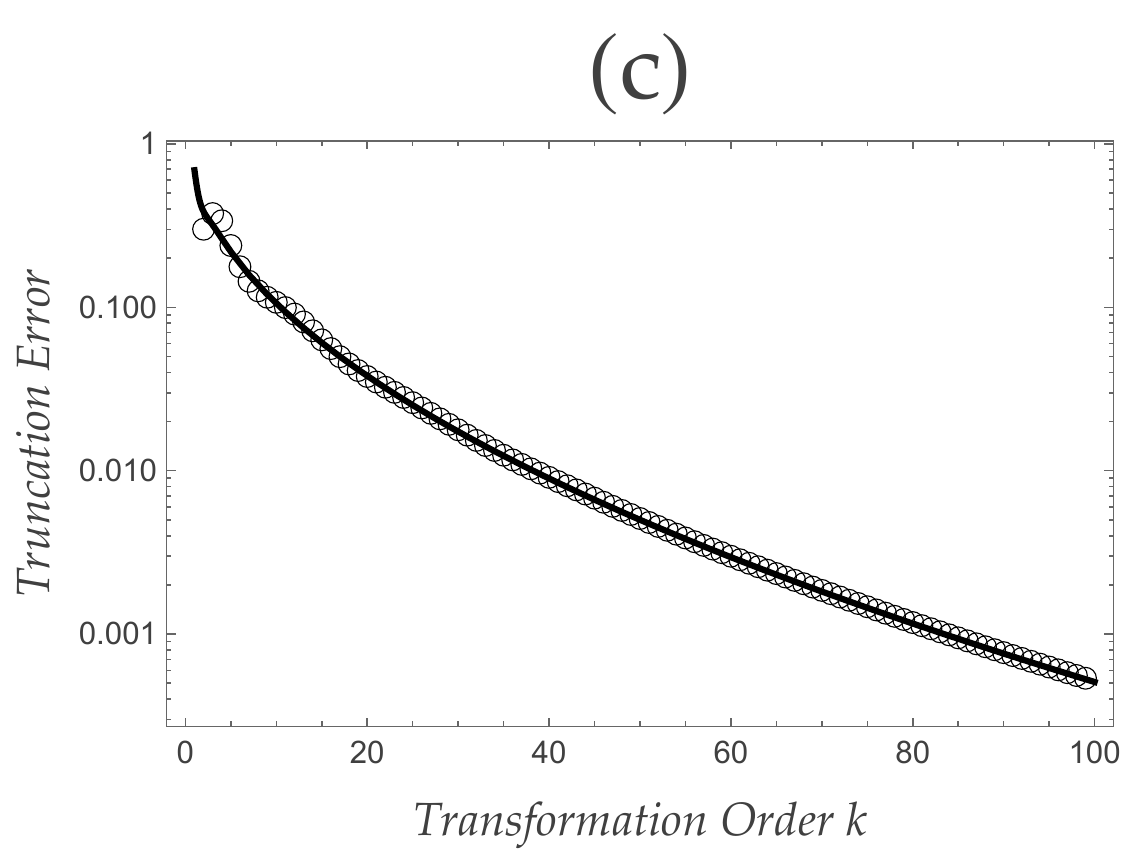}
    \end{minipage}
    \begin{minipage}{6cm}
        \centering
        \includegraphics[width=6cm]{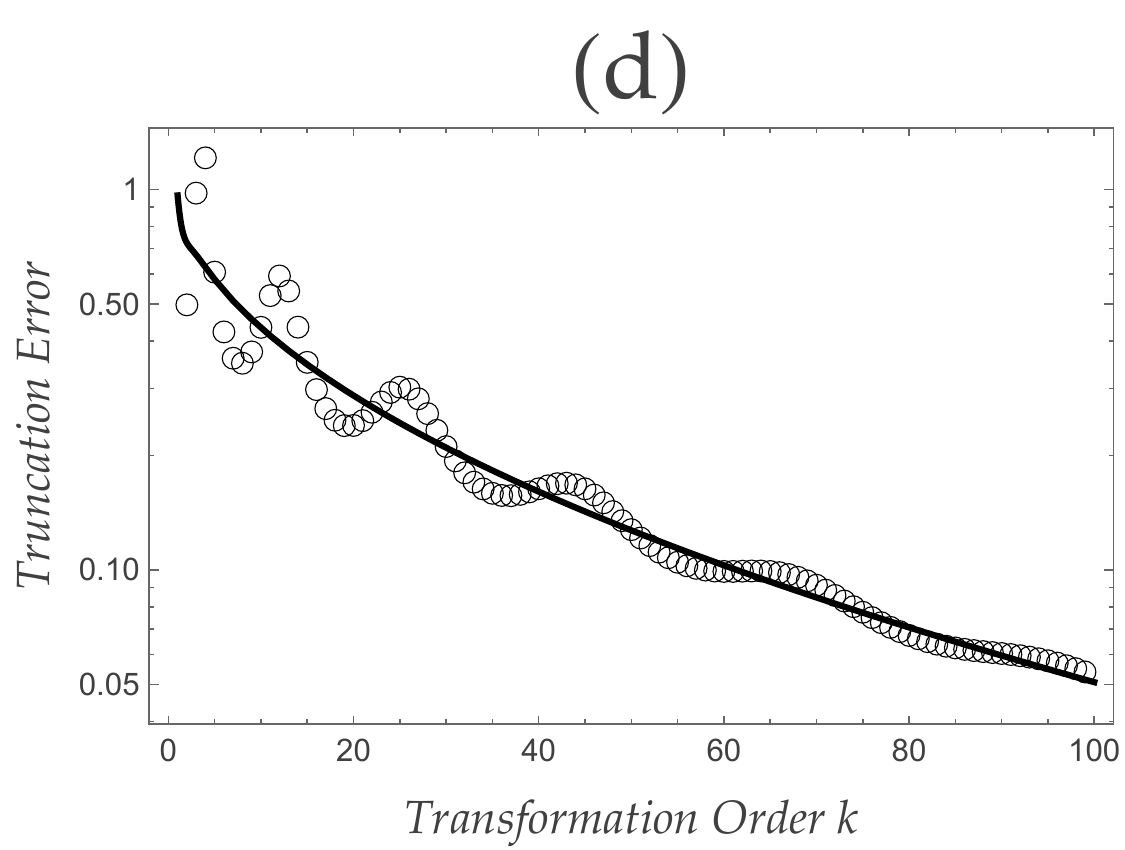}
    \end{minipage}
\caption{The same as in Fig.~\ref{Fig:NumericalResults.1}, but for complex $z=\exp(\mathrm{i}\phi)$,
with $\phi=\pi/4$ (a), $\phi=\pi/2$, (b), $\phi=3\pi/4$ (c), and $\phi=9\pi/10$ (d). }
\label{Fig:NumericalResults.2}
\end{figure}
%


{
A critical benchmark for the proposed method is its performance relative to Padé approximants, 
particularly for series where Carleman's condition is only marginally satisfied. 
}
To this end, we consider  one of the numerical experiments carried out in~\cite{Weniger/Cizek/Vinette/1993}
and already analyzed in Ref.~\cite{Borghi/2025} as far as the sole converging factor was concerned.
On using, for the sake of simplicity, the same notations employed in~\cite{Weniger/Cizek/Vinette/1993,Borghi/2025}, 
consider the following Stieltjes integral:
\begin{equation}
\label{Eq:BorghiAlgorithmImplementation.Target.1.0}
\begin{array}{l}
\displaystyle
\mathcal{J}_3
\,=\,\int_{0}^\infty\,
\dfrac{t^{-1/2}\exp(-t)}{1+\dfrac{64}{45 \pi^2} t^2}\,\mathrm{d} t\,=\,z\,\int_0^\infty\,\dfrac{\mathrm{d}\mu}{z\,+\,t}\,,
\end{array}
\end{equation}
where $z=\dfrac{45 \pi^2}{64}$ and the measure ${\mathrm{d}\mu}$ is given by Eq.~(\ref{Eq:Anharmonic.1.4}) with $q=-1/2$.

In ~\cite{Weniger/Cizek/Vinette/1993}, it was shown numerically how Weniger's $\delta$ transformation (as well as Levin's
$d$) largely outperformed Padé approximants in the resummation of the Stieltjes series which is asymptotic to $\mathcal{J}_3$. 
To help readers, in  Appendix~\ref{App:WynnEpsAlg}, Wynn's $\epsilon$ algorithm is briefly recalled, together with 
its tight relationship with Padé approximants. 
\begin{figure}[ht!]
  \centering
        \includegraphics[width=8cm]{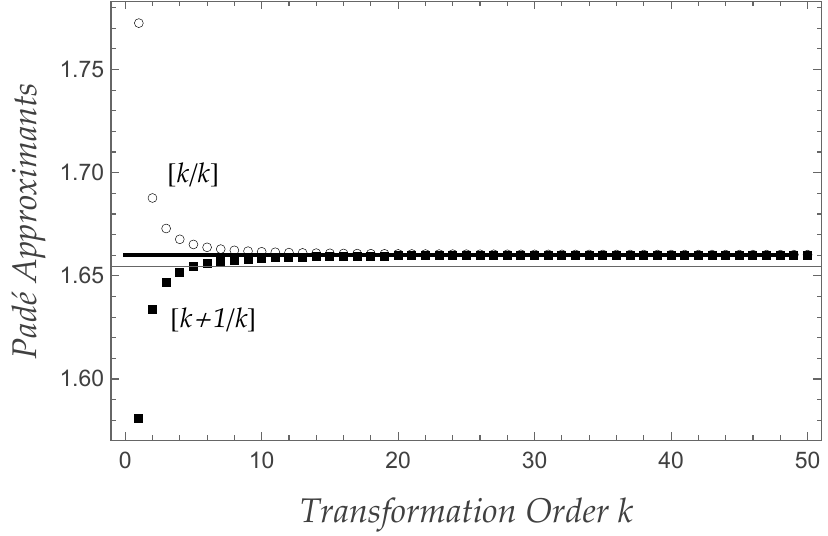}
\caption{Behaviours of Padé approximant sequences $\{[k/k]\}_{k=2}^{50}$ (open circles) and $\{[k+1/k]\}_{k=2}^{50}$ (black squares)
evaluated through Wynn's $\epsilon$-algorithm.
 }
\label{Fig:NumericalResults.3}
\end{figure}

In Fig.~\ref{Fig:NumericalResults.3}, the  behaviors of the Padé approximant sequences $\{[k/k]\}_{k=2}^{50}$ (open circles) and $\{[k+1/k]\}_{k=2}^{50}$ (black squares), obtained through 
Wynn's $\epsilon$-algorithm as $\{\epsilon^{(0)}_{2k}\}_{k=2}^{50}$ and $\{\epsilon^{(1)}_{2k}\}_{k=2}^{50}$, respectively, are plotted against the transformation order $k$: the solid line represents the exact value 
$\mathcal{J}_3 =1.660 177 206 668 046 735 308 9\ldots $
From the figure it is possible to appreciate how the staircase sequence 
$[2/2], [3/2], [3/3], \ldots, [50/50], [51/50]$, provides
lower and upper bounds of $\mathcal{J}_3$. 
More quantitative considerations can be done starting from Eqs. (5.6) - (5.10) of Ref.~\cite{Weniger/Cizek/Vinette/1993}. It follows that, in order to retrieve the Stieltjes integral $\mathcal{J}_3$ with a relative error of the order of $10^{-5}$,  Wynn's algorithm 
needs to be fed by the sequence of the first 100 partial sums of the corresponding asymptotic series. 
The same sequence, when it is $\delta$-transformed, provided a relative error of the order of $10^{-20}$.

{
In order to provide a more quantitative perspective on the convergence patterns, 
in particular to highlight the digits progression towards the exact value of
$\mathcal{J}_3$, 
Table~\ref{tab:comparison} compares the values of the partial sums sequence $\{f_k\}_k$, 
the diagonal Padé approximants sequence $\{[k/k]\}_k$, and the $\delta$ transformation sequence
$\{\delta_k^{(0)}(1/2)\}_k$ for a selected list of transformation order values $k$.
\begin{table}[htpb]
\centering
\caption{Comparison of the numerical approximations of $\mathcal{J}_3$ evaluated at selected orders $k$. The exact value is
rounded to seventeen decimal places,  $\mathcal{J}_3 \simeq 1.66017720666804674$. Digits that agree with the exact value are underlined.}
\label{tab:comparison}
\renewcommand{\arraystretch}{1.3} 
\begin{tabular}{c c c c}
\hline\hline
$k$ & Partial Sums $f_k$ & diagonal Padé $[k/k]$ & $\delta$ transformation $\delta_k^{(0)}(1/2)$ \\
\hline
$10$ & $-4\times 10^7$ & \underline{1.66}144213620852849 & \underline{1.66017}598261438080 \\
$20$ & $-9\times 10^{27}$ & \underline{1.660}59541168291335 & \underline{1.66017720}015209552 \\
$30$ & $- 10^{53}$ & \underline{1.660}39325130457419 & \underline{1.660177206}55817808 \\
$40$ & $-2\times 10^{81}$ & \underline{1.660}31191756303232 & \underline{1.6601772066}7074919 \\
$50$ & $-6\times 10^{111}$ & \underline{1.660}27044462383324 & \underline{1.66017720666}796602 \\
$60$ & $-10^{142}$ & \underline{1.660}24617656962528 & \underline{1.6601772066680}5044 \\
$70$ & $-5\times 10^{177}$ & \underline{1.660}23063190803760 & \underline{1.660177206668046}51 \\
$80$ & $- 5 \times 10^{212}$ & \underline{1.660}22001562972316 & \underline{1.6601772066680467}5 \\
$90$ & $-6\times 10^{248}$ & \underline{1.660}21240959133277 & \underline{1.6601772066680467}3 \\
$100$ & $-7\times 10^{285}$ & \underline{1.660}20675400358820 & \underline{1.66017720666804674} \\
\hline\hline
\end{tabular}
\end{table}

As Table~\ref{tab:comparison} clearly illustrates, while the sequence of partial sums rapidly explodes, Padé approximants slowly lock onto the exact value, struggling to achieve high precision. On the opposite, 
the $\delta$ transformation consistently captures a significantly higher number of correct digits, confirming its computational superiority for this class of Stieltjes series.
}

{
It had already been shown  in~\cite{Borghi/Weniger/2015} that the Padé approximant transformation error for the Euler series 
turns out to be asymptotically proportional to $\exp(-4 z^{1/2} k^{1/2})$, in agreement with \cite[Eq.~(8.5.18b)]{Bender/Orszag/1978}. This result was found on using the fact that Euler series' diagonal Padé approximants can be computed via Drummond's transformation~\cite{Sidi/1981}.
Unfortunately, for the Stieltjes series under investigation, no similar analytical results have been found. 
For this reason, the comparison of the resummation performances of the $\delta$ transformation with those of 
Pad\'e approximants will be carried out numerically. 
}
\begin{figure}[ht!]
  \centering
        \includegraphics[width=8cm]{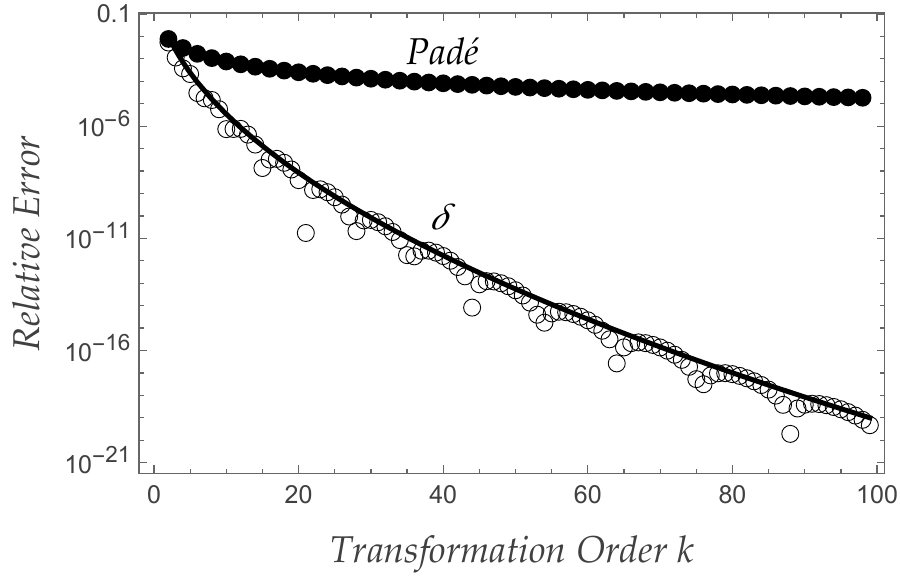}
\caption{
Behaviour of the relative error, against the transformation order $k$, provided 
by the sequence $\{\delta_k^{(0)}(1/2)\}_{k=2}^{100}$ (open circles) and by the diagonal sequence of Padé approximants $\{[k/k]\}_{k=2}^{100}$ (open squares),  for the resummation of the Stieltjes function corresponding to $(z,q)=\left(\dfrac{45}{64}\pi^2,-\dfrac 12\right)$. 
The solid curve represents the asymptotic estimate $\exp(-2\sqrt 2 \,z^{1/4}\, k^{1/2})$.
}
\label{Fig:NumericalResults.4}
\end{figure}

In Fig.~\ref{Fig:NumericalResults.4}, the relative errors provided by the $\delta$ transformation (open circles) and by diagonal Padé approximants (open squares) are shown against the transformation order $k$, for $z=45\pi^2/64$. Differently from Fig.~\ref{Fig:NumericalResults.1}, the sinusoidal factor in the  analytical estimate of Eq.~(\ref{Eq:WenigerTransformationErrorAsymptotics})  has been omitted, in order to emphasize the asymptotic dominance of the exponential factor $\exp \left(-2\sqrt 2\, z^{1/4} \, k^{1/2}\right)$, which is independent of the value of $q$.

{
Figure~\ref{Fig:NumericalResults.4} highlights the different convergence regimes of the two methods. 
While the $\delta$ transformation exhibits the predicted exponential decay, the convergence of Padé approximants is severely hindered. As previously noted, achieving high-precision resummation with Padé methods for this specific class of series is computationally demanding. Numerical data in Fig.~\ref{Fig:NumericalResults.4} indicates a  subexponential convergence rate; numerical fits of the relative error using functions of the form $A \exp(-B k^m)$ proved unsuccessful.
}

{
While a full theoretical investigation of the Padé convergence rate for this specific problem is beyond the scope of the present paper, this behavior can be mathematically framed by examining Carleman's condition. As is well known, the divergence of Carleman's series in Eq.~(\ref{Eq:Calerman}) can be easily proved by using Stirling's formula, i.e.,
}
{%
\begin{equation}
\label{Eq:CalermanAsymptotic}
\displaystyle
\mu_m = \Gamma (2m+q+1)\,\sim\left(2m\right)!\,\sim\,\sqrt{4\pi\,m}\,\left(\dfrac{2m}{\rm e}\right)^{2m},
\quad m \to \infty\,,
\end{equation}
so that
\begin{equation}
\label{Eq:CalermanAsymptotic.2}
\displaystyle
\mu_m^{-\frac 1{2m}} \,\sim\, 
(4\pi\,m)^{-\frac 1{4m}}\,\left(\dfrac{2m}{e}\right)^{-1} \,\sim\, 
\dfrac e{2m},
\qquad\qquad m \to \infty.
\end{equation}
This, in turn,  proves the validity of Carleman's condition, since the series 
$\displaystyle\sum_m\,\mu_m^{-\frac1{2m}}$ {does diverge}, although logarithmically,
i.e., like the harmonic series.  
In the case of Euler' series ($\mu_m=m!$), the Carleman series
diverges as $\displaystyle\sum_m m^{-1/2}$, a strongest divergence condition with respect to 
the harmonic one.
More generally speaking, if $\mu_m=((2+\epsilon)m)!$, the 
Carleman series asymptotically behaves like $\displaystyle\sum_m \frac 1{m^{1+\epsilon/2}}$, i.e.,
it converges for $\epsilon >0$ and diverges for $\epsilon \le 0$. In other words, our Stieltjes series 
is placed at the edge between divergence and convergence of the associated Carleman series.

Another example of Stieltjes series whose Carleman series displays a logarithmical divergence is characterized by the 
moment sequence $\mu_m=(m!)^2$~\cite[Corollary 12.11h ]{Henrici/1977}. 
Accordingly, it could be worth comparing the action of $\delta$ and Padé in this case, similarly as it was done in Fig.~\ref{Fig:NumericalResults.4}.  In the present case 
the associated measure $\mathrm{d}\mu$ turns out to be 
\begin{equation}
\label{Eq:CalermanAsymptotic.2.0.1}
\displaystyle
\mathrm{d}\mu\,=\,2 K_0(2\sqrt t)\,\mathrm{d}t,\qquad\qquad t\in [0,\infty),
\end{equation}
{
where $K_0(\cdot)$ denotes the modified Bessel function of the second kind. 
The corresponding Stieltjes integral  can expressed as follows:
\begin{equation}
\label{Eq:CalermanAsymptotic.2.1}
\displaystyle
\int_0^\infty\,\dfrac{\mathrm{d}\mu}{z+t}\,=\,
\int_0^\infty\, \dfrac{2K_0(2\sqrt{t})}{z+t}\,\mathrm{d}t\,=\,
4\,\int_0^\infty\, \dfrac{u}{u^2+z}\,K_0(2u)\,\mathrm{d}u,
\end{equation}
where the new integration variable $u=\sqrt t$ has been introduced. The integral in Eq.~(\ref{Eq:CalermanAsymptotic.2.1})
can be solved via~\cite[formula 2.16.3.14]{Prudnikov/Brychkov/Marichev/1986b}, which gives
\begin{equation}
\label{Eq:App:K0.2}
\begin{array}{l}
\displaystyle
\int_0^\infty\,\dfrac{\mathrm{d}\mu}{z+t}\,=\,
4\,\mathcal{S}_{-1,0}(2\sqrt z),
\end{array}
\end{equation}
where the symbol $\mathcal{S}_{\alpha,\beta}(\cdot)$ denotes the so-called Lommel function of the second kind~\cite{DLMF/Book/2010}.
}

Figure~\ref{Fig:NumericalResults.5} shows the same numerical experiment carried out in Fig.~\ref{Fig:NumericalResults.4},
but for the Stieltjes series associated to the measure in Eq.~(\ref{Eq:CalermanAsymptotic.2.0.1}). From this figure 
{
a clear similarity in the relative error behaviors can be appreciated: the $\delta$ transformation maintains its strong exponential convergence, while Padé approximants suffer from the same severe slowdown observed in the $(2n)!$ case. This confirms the suitability and the computational robustness of the $\delta$ transformation for decoding Stieltjes series situated exactly at the boundary of Carleman's condition.
}
\begin{figure}[ht!]
  \centering
        \includegraphics[width=8cm]{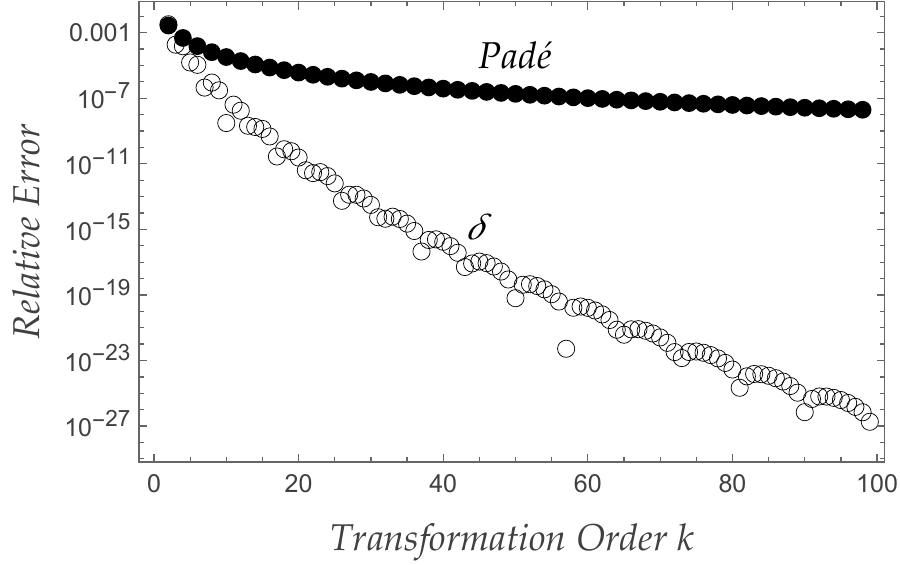}
\caption{
Behaviour of the relative error, against the transformation order $k$, provided 
by the sequence $\{\delta_k^{(0)}(1)\}_{k=2}^{100}$ (open circles) and by the diagonal sequence of Padé approximants $\{[k/k]\}_{k=2}^{100}$ (open squares),  for the resummation of the Stieltjes function in Eq.~(\ref{Eq:CalermanAsymptotic.2.1}) when $z=\dfrac{45}{64}\pi^2$. 
}
\label{Fig:NumericalResults.5}
\end{figure}
}

\section{Conclusions}
\label{Sec:Conclusions}

{
In the present paper, we have provided an asymptotic analysis of Weniger’s $\delta$ transformation applied to a class of superfactorially divergent Stieltjes series, whose moments grow as $(2n)!$. By using some recent results on converging factors, 
an integral representation for the truncation error has been derived, which allowed us to derive a leading-order asymptotic 
estimate of the related convergence rate. 

These theoretical findings have extensively been validated through several numerical experiments. Our results demonstrate that the 
$\delta$ transformation significantly outperforms Padé approximants for series situated at the boundary of Carleman’s condition, where traditional resummation methods exhibit a transition from exponential to sub-exponential convergence. 
Moreover, the extension of our asymptotic formula to the complex plane, presented here only as an analytical conjecture, 
has been supported by numerical evidence, showing that the transformation remains robust and efficient deep into the left 
half-plane, provided the Stieltjes function branch cut is not crossed.

Our theoretical development addresses a significant gap still present in the literature, where the mathematical understanding 
of nonlinear Levin-type transformations has remained largely underdeveloped when compared to classical Padé theory. 
Weniger's $\delta$ transformation proves to be a powerful and computationally efficient tool for the resummation of $(2n)!$-divergent expansions. The mathematical framework developed in this work provides a deeper understanding of the convergence mechanism of this particular type of nonlinear sequence transformation.

Nevertheless, it should be emphasized how the transition from factorial to superfactorial divergence introduces 
substantial analytical hurdles, which made the technical challenges tackled in the present paper considerably harder 
than those addressed in Ref.~\cite{Borghi/Weniger/2015}. The present work represents only the second contribution 
to a remarkably sparse list of analytical results for $\delta$ transformation.

While real-world physical applications often lead to nonStieltjes series with important complex subdominant structures
that can elude convergence, a point raised in Barry Simon's 1971 critique of Padé theory~\cite{Simon/1971}, we believe that 
establishing a rigorous framework within the Stieltjes domain for Levin-type transformations represents a fundamental and 
necessary theoretical milestone upon which future generalizations can be built.
In this respect, the analytical results established here provide, for the superfactorial Stieltjes class here explored, a bridge between 
the numerical efficiency of Weniger’s transformation and the rigorous theoretical foundation required for its reliable 
application in physics and applied mathematics. We hoped that the present findings will stimulate further investigations 
into the resummation of even more aggressively divergent series, also including nonStieltjes scenarios.
}

\acknowledgments{I wish to thank Turi Maria Spinozzi for his useful comments.}

\appendix

\section[\appendixname~\thesection]{Proof of Eq.~(\ref{Eq:NumeratoreIntegrale})}
\label{App:Eq:NumeratoreIntegrale}
{
We start by substituting from Eq.~(\ref{FS_IntRep}) into the numerator of Eq.~(\ref{Eq:WenigerTransformationOnK.2}).
On replacing $n$ by $n+1$, the integral representation in Eq.~(\ref{FS_IntRep}) yields $\varphi_{n+1}=\int_0^1\,t^{n+\beta} \Phi(t)\mathrm{d}t$. In order to simplify the subsequent notations within the finite difference operators, we set $\beta+1=\gamma$,
so that $t^{n+\beta}=t^{n+\gamma-1}$. This gives 
}
\begin{equation}
  \label{Eq:App:Eq:NumeratoreIntegrale.1}
\begin{array}{l}
\displaystyle
\Delta^k\left\{(n+\gamma)_{k-1}\varphi_{n+1}\right\}\,=\,
\Delta^k\left\{(n+\gamma)_{k-1}\int_{0}^{1} \, t^{n+\gamma-1} \, \Phi (t) \, \mathrm{d} t \right\} \,=\,\\
\\
\,=\,
\displaystyle
\int_{0}^{1}\, \mathrm{d} t \, t^{\gamma-1}\,\Phi (t) \,\Delta^k\left\{(n+\gamma)_{k-1} t^n \right\},
\end{array}
 \end{equation}
Now, we have~\cite{Borghi/Weniger/2015}
\begin{equation}
  \label{Eq:App:Eq:NumeratoreIntegrale.2}
\begin{array}{l}
\displaystyle
\Delta^k\left\{(n+\gamma)_{k-1} t^n \right\}\,=\,
(-1)^k\,
(n+\gamma)_{k-1} \, t^{n} \, 
{}_{2} F_{1}
\left(
    \genfrac{}{}{0pt}{}{-k,k+n+\gamma-1} {n+\gamma}; t \right) \, ,
\end{array}
 \end{equation}
where ${}_2F_1(\cdot)$ denotes the Gauss hypergeometric function. 
Finally, on substituting from Eq.~(\ref{Eq:App:Eq:NumeratoreIntegrale.2}) into Eq.~(\ref{Eq:App:Eq:NumeratoreIntegrale.1}), 
after simple algebra Eq.~(\ref{Eq:NumeratoreIntegrale}) follows.

\section[\appendixname~\thesection]{Analytical evaluation of the Stieltjes integral in Eq.~(\ref{Eq:Anharmonic.1.4})}
\label{App:AnalyticalStieltjes}

{
On performing the  variable substitution $u=\sqrt t$ into the integral into Eq.~(\ref{Eq:Anharmonic.1.4}) we have 
\begin{equation}
\label{Eq:App:AnalyticalStieltjes.1}
f(z)\,=\,
\int_0^\infty\,\dfrac{u^{q}\,\exp(-u)}{u^2+z}\,\mathrm{d}u.
\end{equation}
The integral can be evaluated, for $q >-1$ and $z>0$, as follows. First of all, 
the integrand is recast as
\begin{equation}
\label{Eq:App:AnalyticalStieltjes.2}
f(z)\,=\,
\dfrac {\mathrm{i}}{2\sqrt{z}}
\left[
\int_0^\infty\,\dfrac{u^{q}\,\exp(-u)}{u+\mathrm{i}\sqrt z}\,\mathrm{d}u\,-\,
\int_0^\infty\,\dfrac{u^{q}\,\exp(-u)}{u-\mathrm{i}\sqrt z}\,\mathrm{d}u
\right],
\end{equation}
then, on using~\cite[formula 2.3.6.13]{Prudnikov/Brychkov/Marichev/1986a}, simple algebra  gives 
\begin{equation}
\label{Eq:App:AnalyticalStieltjes.3}
f(z)\,=\,
-\dfrac {\Gamma(q+1)}{\sqrt{z}}
\mathrm{Im}\left\{
(-z)^{q/2}\,\exp(\mathrm{i}\sqrt z)\,\Gamma(-q,\mathrm{i}\sqrt z)
\right\}, \qquad q>-1,\,\,z>0.
\end{equation}
}

\section[\appendixname~\thesection]{Proof of Eq.~(\ref{Eq:Thm:AsymptoticsDenominatorDelta.1})}
\label{App:Eq:Thm:AsymptoticsDenominatorDelta.1}

Our proof starts from the following expansion~\citep[Eq.\ 7.4.5(6) on p.\ 263]{Luke/1969a}:
\begin{align}
  \label{Luke/1969a_Eq_7.4.5(6)}
   & {}_{p+2} F_{q} \left(
    \genfrac{}{}{0pt}{}{-n, n+\lambda, \alpha_{p}} 
     {\rho_{q}}; -z \right) \; \sim \; \sum_{t=1}^{p} \,
      \frac{(n+\lambda)_{-\alpha_{t}}}{(n+1)_{\alpha_{t}}} \,
       \mathcal{L}_{p+2, q}^{(\alpha_{t})} 
        (z \mathrm{e}^{\mathrm{i} \delta \pi})
  \notag
  \\
  & \qquad + \frac {(2\pi)^{(1-\beta)/2} \Gamma (\rho_{q})}{\beta^{1/2}
    \Gamma (\alpha_{p})} \, \bigl( N^{\beta} z \bigr)^{\nu} \,  
     \exp \left( N z^{1/\beta} \beta - (a z/3) - \Omega (-z)/
      (N z^{1/\beta}) + \mathrm{O} (N^{-2}) \right) \, ,
  \notag
  \\
  & \qquad \qquad \vert \arg (z) \vert \le \pi - \epsilon \, , 
   \quad \epsilon > 0 \, , \quad  \delta = +(-) \quad \text{if} \quad
    \arg (z) \le (>) \; 0 \, .
\end{align}
In order to apply Eq.~(\ref{Luke/1969a_Eq_7.4.5(6)}) to the present case, we have to extract the right values of the various parameters. In this way we obtain $n=k$, $p=0$, $q=3$, $\lambda=q/2$, and $\{\rho_q\}=\left\{\dfrac q2+1,\dfrac{q}{2}+\dfrac{3}{2},\dfrac{q}{2}+2\right\}$. Accordingly, the finite sum is not present (since $p=0$).
The other parameters to be inserted into Eq.~(\ref{Luke/1969a_Eq_7.4.5(6)}) can be derived directly from the prescriptions given in~\citep[Chapter 7.4]{Luke/1969a}. More precisely, we have $\beta=4$ and $a=0$, while
\begin{equation}
\begin{array}{l}
N^4=k\left(k+\dfrac q2\right)\,\Longrightarrow N=\left[k\left(k+\dfrac q2\right)\right]^{1/4}\,.
\end{array}
 \end{equation}
Moreover, 
\begin{equation}
\left\{
\begin{array}{l}
\Gamma(\rho_q)\,=\,\Gamma\left(\dfrac q2+1\right)\,\Gamma\left(\dfrac q2+\dfrac 32\right)\,\Gamma\left(\dfrac q2+2\right),\\
B_1=0, 
\\
C_1=\displaystyle\sum_{t=1}^3\rho_t=\dfrac q2+1+\dfrac q2+\dfrac 32+\dfrac q2+2=\dfrac 32 (3+q),
\\
\gamma=\dfrac {\beta-1+2B_1-2C_1}{2\beta}=-\dfrac 38(2+q),
\end{array}
\right.
\end{equation}
so that Eq.~(\ref{Luke/1969a_Eq_7.4.5(6)}) becomes
\begin{align}
  \label{Luke/1969a_Eq_7.4.5(6)ParticularCase}
   & {}_{2} F_{3} \left(
    \genfrac{}{}{0pt}{}{-k, k+\dfrac q2} 
     {\rho_{q}}; -\dfrac z4 \right) \; \sim \;\frac {(2\pi)^{(1-4)/2} 
\Gamma\left(\dfrac q2+1\right)\,\Gamma\left(\dfrac q2+\dfrac 32\right)\,\Gamma\left(\dfrac q2+2\right)}
{4^{1/2}}
  \notag
  \\
  & \times \left[k\left(k+\dfrac q2\right)\dfrac z4\right]^{-\frac 38(q+2)}\,  
\exp \left\{4 \left[k\left(k+\dfrac q2\right) \dfrac z4\right]^{1/4}+ \mathcal{O}(N^{-1}) \right\} \, ,
  \end{align}
If we now discard, in the exponential factor in Eq.\ (\ref{Luke/1969a_Eq_7.4.5(6)ParticularCase}), 
all contributions that vanish at least like $\mathcal{O} \bigl( 1/N \bigr)$ as $N \to \infty$, after
simple algebra the following leading order asymptotic approximation for the denominator in 
Eq.~(\ref{Eq:Thm:AsymptoticsDenominatorDelta.1}) is obtained:
\begin{equation}
  \label{Eq:TransfErrorIntegraleN0Beta1Denominator}
\begin{array}{l}
\displaystyle
{}_{2} F_{3} \left(
    \genfrac{}{}{0pt}{}{-k, k+\dfrac q2} 
     {\rho_{q}}; -\dfrac z4 \right) \; \sim \;
\\
\\
\sim
\dfrac { 
\Gamma\left(\dfrac q2+1\right)\,\Gamma\left(\dfrac q2+\dfrac 32\right)\,\Gamma\left(\dfrac q2+2\right)}
{2(2\pi)^{3/2}}\,\left[k\left(k+\dfrac q2\right)\dfrac z4\right]^{-\frac 38(q+2)}\,\exp \left\{4 \left[k\left(k+\dfrac q2\right) \dfrac z4\right]^{1/4}\right\}  \, ,\qquad \qquad\qquad \qquad\qquad \qquad\qquad \qquad  k \to \infty
\end{array}
 \end{equation}
{ or, since $k \gg 1 >|q|$,} 
\begin{equation}
  \label{Eq:TransfErrorIntegraleN0Beta1DenominatorBis}
\begin{array}{l}
\displaystyle
{{}_{2} F_{3} \left(
    \genfrac{}{}{0pt}{}{-k, k+\dfrac q2} 
     {\rho_{q}}; -\dfrac z4 \right) \; \sim \;}\\
\\
\sim
\dfrac { 
\Gamma\left(\dfrac q2+1\right)\,\Gamma\left(\dfrac q2+\dfrac 32\right)\,\Gamma\left(\dfrac q2+2\right)}
{2(2\pi)^{3/2}}\,\left(\dfrac z4\right)^{-\frac 38(q+2)}\,k^{-\frac 34(q+2)}\, \exp \left\{2\sqrt 2 k^{1/2} z^{1/4}\right\}  \, ,\qquad \qquad\qquad \qquad\qquad \qquad\qquad \qquad  k \to \infty
\end{array}
 \end{equation}
which coincides with Eq.~(\ref{Eq:Thm:AsymptoticsDenominatorDelta.1}).

\section[\appendixname~\thesection]{Proof of Eq.~(\ref{Eq:Thm:AsymptoticsNumeratorDelta.1})}
\label{App:Eq:Thm:AsymptoticsNumeratorDelta.1}

First of all, we have to find the asymptotics of the hypergeometric function ${}_2F_1$. To this end, it is worth
invoking~\cite[Eq. 7.4.2(8) on p.250]{Luke/1969a}, which gives
\begin{align}
  \label{Luke/1969a_Eq_7.4.2(8)}
  & {}_{p+2} F_{p+1} \left(
   \genfrac{}{}{0pt}{}{-n, n+\lambda, \alpha_{p}}      
    {\rho_{p+1}}; z \right) \; \sim \; 
    \sum_{t=1}^{p} \, \frac 
     {(n+\lambda)_{-\alpha_{t}}}{(n+1)_{\alpha_{t}}} \,
      \mathcal{L}_{p+2, p+1}^{(\alpha_{t})} (z) \, + \, 
       \frac {\Gamma (\rho_{p+1}) \, N^{2\gamma}} 
        {\Gamma (\alpha_{p}) \, \Gamma(1/2)} \, 
         \frac {[\sin (\theta/2)]^{2\gamma}}
          {[\cos (\theta/2)]^{2\gamma+\lambda}}    
  \notag
  \\
  & \quad \times 
   \exp \left\{ \frac{\varphi_{2} (\theta) + a_{2}}{N^{2}} + 
    \mathrm{O} \left( \frac{1}{N^{4}} \right) \right\} \, 
     \cos \left\{ N \theta + \pi \gamma + \frac{\varphi_{1} (\theta)}{N} 
       \, + \, \frac{\varphi_{3} (\theta)} {N^{3}} \, + \, 
        \mathrm{O} \left( \frac{1}{N^{5}} \right) \right\} \, ,
  \notag
  \\
  & \qquad \vert \arg (z) \vert \le \pi - \epsilon \, , 
   \quad \vert \arg (1-z) \vert \le \pi - \epsilon \, , 
	\quad \epsilon \; > \; 0 \, .
\end{align}
Again, we shall set $n=k$ and $\lambda=q/2$, while $N =
\sqrt{k(k+q/2)}$ and $z=t \in [0,1]$. In addition, we have $p=0$, so that, as in Eq.\ (\ref{Luke/1969a_Eq_7.4.5(6)}), the finite
sum  disappears.
Moreover, on neglecting all contributions that vanish at least like $\mathcal{O} \bigl( 1/N \bigr)$ as
$N \to \infty$, we have that the exponential factor 
can be approximated by one, while the cosinusoidal factor reduces to $\cos(N\theta+\pi \gamma)$,
where $\theta$ is related to $t$ by $t=\sin^{2} (\theta/2)$, so that $\theta\in [0,\pi]$. The remaining unspecified
quantities are defined in~\citep[Eq.\ 7.4.2(9) on pages\ 251 and\
252.]{Luke/1969a}. More precisely, we have
\begin{equation}
\left\{
\begin{array}{l}
B_1=0, 
\\
C_1=1+\dfrac  q2,
\\
\gamma=\dfrac {1+2B_1-2C_1}{4}=-\dfrac {1+q}4,
\end{array}
\right.
\end{equation}
so that
\begin{align}
  \label{Eq:TransfErrorIntegraleN0Beta1NumeratorBis}
  & {}_{2} F_{1} \left(
   \genfrac{}{}{0pt}{}{-k, k+\dfrac q2}      
    {\rho_{1}}; t \right) \; \sim \; 
       \dfrac {\Gamma (\rho_{1}) \, N^{2\gamma}} 
        {\Gamma(1/2)} \, 
         \frac {[\sin (\theta/2)]^{2\gamma}}
          {[\cos (\theta/2)]^{2\gamma+\lambda}}    \cos \left\{ N \theta + \pi \gamma \right\}\,,\qquad k\to\infty,
  \notag
\end{align}
i.e.,
\begin{equation}
  \label{Eq:TransfErrorIntegraleN0Beta1NumeratorBis.1}
\begin{array}{l}  
{}_{2} F_{1} \left(
   \genfrac{}{}{0pt}{}{-k, k+\dfrac q2}      
    {\rho_{3}}; t \right) \; \sim \;\\
\\
\sim
\displaystyle
\dfrac { 
\Gamma\left(\dfrac q2+1\right)}
{\sqrt\pi} 
\left[k(k+q/2)]\right]^{-\frac {1+q}4}\,
\displaystyle
\quad \times 
         \frac {[\sin (\theta/2)]^{2\gamma}}
          {[\cos (\theta/2)]^{2\gamma+\lambda}} \cos \left\{ N \theta + \pi \gamma \right\}\,,
\qquad k\to\infty,
\end{array}
\end{equation}
or, on taking into account that $2\gamma+\lambda= -\dfrac{1+q}2+\dfrac q2=-\dfrac 12$,
\begin{equation}
  \label{Eq:TransfErrorIntegraleN0Beta1NumeratorBis.2}
\begin{array}{l} 
 {}_{2} F_{1} \left(
   \genfrac{}{}{0pt}{}{-k, k+\dfrac q2}      
    {1+\dfrac q2}; t \right) \; \sim \;
\dfrac {\Gamma\left(\dfrac q2+1\right)}{\sqrt\pi} 
\,k^{-\frac {1+q}2}\,
         \sin^{-\frac{1+q}2} \dfrac\theta 2\,\sqrt{\cos \dfrac\theta 2}
          \cos \left\{ k \theta - \dfrac{1+q}4\pi\right\}\,,
\qquad k\to\infty,
\end{array}
\end{equation}
which coincides with Eq.~(\ref{Eq:Thm:AsymptoticsNumeratorDelta.1}).

\section[\appendixname~\thesection]{Derivation of Eq.~(\ref{Eq:ExactNumeratorAsymtpt.1})}
\label{App:ExactNumeratorAsymtpt.1}

We start on recasting the integrand of Eq.~(\ref{Eq:ExactNumeratorAsymtpt.1}) as the product $\Psi(\theta)\,\Lambda(\theta)$, where
\begin{equation}
  \label{App:Eq:ExactNumeratorAsymtpt.Formula.1}
\begin{array}{l}
\displaystyle
\Psi(\theta)\,=\,
\sin\left[\sqrt z\,\left(\dfrac 1{\sin\dfrac\theta 2}\,-\,1\right)\right]\,
\cos\left(k\theta-\eta\dfrac \pi2\right),
\end{array}
 \end{equation}
and
\begin{equation}
  \label{App:Eq:ExactNumeratorAsymtpt.Formula.2}
\begin{array}{l}
\displaystyle
A(\theta)\,=\,
\sin^{\eta}\dfrac\theta 2\,\cos^{3/2}\dfrac\theta 2,
\end{array}
 \end{equation}
where for the sake of convenience the parameter $\eta=\dfrac{1+q}2 \in [0,1]$, has been introduced. 
On applying Werner formulas, $\Psi$ can be rewritten as follows:
\begin{equation}
  \label{App:Eq:ExactNumeratorAsymtpt.Formula.3}
\begin{array}{l}
\displaystyle
\Psi(\theta)\,=\,\dfrac 12\,\Psi_+(\theta)-\dfrac 12\,\Psi_-(\theta),
\end{array}
 \end{equation}
where
\begin{equation}
  \label{App:Eq:ExactNumeratorAsymtpt.Formula.4}
\begin{array}{l}
\displaystyle
\Psi_\pm(\theta)\,=\,
\sin\left[
\phi_\pm(\theta)
-\eta\dfrac \pi2
\right],
\end{array}
 \end{equation}
with
\begin{equation}
  \label{App:Eq:ExactNumeratorAsymtpt.Formula.4.1}
\begin{array}{l}
\displaystyle
\phi_\pm(\theta)\,=\,
k\theta\pm\sqrt z\left(\dfrac 1{\sin\dfrac\theta 2}\,-\,1\right)\,.
\end{array}
 \end{equation}
Accordingly, the starting integral, say $\mathcal{I}$, can be written as 
$\mathcal{I}=\mathcal{I}_+-\mathcal{I}_-$, where
\begin{equation}
  \label{App:Eq:ExactNumeratorAsymtpt.Formula.5}
\begin{array}{l}
\displaystyle
\mathcal{I}_\pm\,=\,\dfrac 12\,\int_0^\pi\,A(\theta)\,\Psi_\pm(\theta)\,\mathrm{d}\theta\,.
\end{array}
 \end{equation}
The idea consists in estimating (if possible) the two integrals on using standard stationary phase techniques. In particular, the factor $\Psi_\pm(\theta)$ represents the highly oscillating function, due to the presence, for $k \to \infty$, of $k\theta$ as well as to the presence of the term
$1/\sin\frac\theta 2$ (in the limit of $\theta\to 0$). The function $A(\theta)$
contains the slowly varying factor. Accordingly, to  asymptotically estimate $\mathcal{I}_\pm$,
the real solutions of the following equations have to be found:
\begin{equation}
  \label{App:Eq:ExactNumeratorAsymtpt.Formula.6}
\begin{array}{l}
\displaystyle
\dfrac{\mathrm{d}}{\mathrm{d}\theta}	\,\phi_\pm(\theta)\,=\,0,
\end{array}
 \end{equation}
i.e., on taking Eq.~(\ref{App:Eq:ExactNumeratorAsymtpt.Formula.4.1}) into account,
\begin{equation}
  \label{App:Eq:ExactNumeratorAsymtpt.Formula.6.1}
\begin{array}{l}
\displaystyle
k\,=\,\pm\,\dfrac{\sqrt z}2\,\dfrac{\cos\dfrac\theta 2}{\sin^2\dfrac\theta 2}\,.
\end{array}
 \end{equation}
{
The two contributions will now be treated separately. 

As far as $\mathcal{I}_+$ is concerned, it is characterized by the presence of a stationary 
point within the integration domain. The equation $\phi'_+(\theta_s)=0$ gives the position of 
a saddle point which, for large values of $k$, is located approximately at 
\begin{equation}
  \label{App:Eq:ExactNumeratorAsymtpt.Formula.7.1}
\begin{array}{l}
\displaystyle
\theta_s \sim \dfrac{\sqrt 2\,z^{1/4}}{\sqrt k}. 
\end{array}
 \end{equation}
For $k \to \infty$, the stationary point approaches the lower endpoint of the integration domain $\theta=0$, which 
could imply a potential coalescence between the saddle and the endpoint with a consequence need for uniform expansions.
To understand why stationary phase approximation keeps its validity in the present scenario, it is sufficient to recast 
$\mathcal{I}_+$ in terms of the new integration variable $\xi=\sqrt k\,\theta$, which maps the integration domain from 
$\theta \in [0,1]$ into $\xi \in [0,\sqrt k]$. On extending, in the asymptotic limit $k\to\infty$, the $\xi$-integration domain to the semi-infinite interval $[0,\infty)$, we have
\begin{equation}
  \label{App:Eq:ScaledIntegral}
\begin{array}{l}
\displaystyle
\mathcal{I}_+\,\sim\,\dfrac 1{\sqrt{k}}\,
\int_0^\infty\,\mathrm{d}\xi\,
A\left(\dfrac\xi{\sqrt k}\right)\,
\sin\left[\phi_+\left(\dfrac\xi{\sqrt k}\right) -\eta\dfrac \pi 2\right],
\qquad\qquad k \to\infty,
\end{array}
 \end{equation}
where the scaled phase now becomes, in the limit $k\to\infty$,
\begin{equation}
  \label{App:Eq:ScaledIntegral.2}
\begin{array}{l}
\displaystyle
\phi_+\left(\dfrac\xi{\sqrt k}\right)\,=\,
k\dfrac \xi{\sqrt k}+\sqrt z\left(\dfrac 1{\sin\dfrac \xi{2\sqrt k}}\,-\,1\right)\,\sim\,
\sqrt k\, \xi+\sqrt z\left(\dfrac {2\sqrt k}{\xi}\,-\,1\right)\,\sim\,
\sqrt k \left(\xi+\dfrac{2\sqrt z}\xi\right). 
\end{array}
 \end{equation}
%
The stationary point of the integral in the $\xi$-domain
$[0,\infty]$ now turns out to be $\xi_s=\sqrt 2\,z^{1/4}$, 
which implies that the distance between the saddle and the boundary remains constant with the saddle well separated from the origin. Stationary phase method
can safely  be applied to extract the leading order of $\mathcal{I}_+$. To this end, it is worth coming back into 
the $\theta$-domain (the math is simpler), where 
\begin{equation}
\label{App:Eq:ExactNumeratorAsymtpt.Formula.8}
\begin{array}{l}
\displaystyle
\mathcal{I}_+\,\sim\,
\dfrac {A(\theta_s)}2\,\sqrt{\dfrac{2\pi}{|\phi''(\theta_s)|}}\,
\sin\left[\phi(\theta_s)+\dfrac\pi 4\,\mathrm{sign}\{\phi''(\theta_s)\}\,-\,\eta\dfrac\pi 2\right].
\end{array}
 \end{equation}
From Eq.~(\ref{App:Eq:ExactNumeratorAsymtpt.Formula.4.1}) we have
\begin{equation}
  \label{App:Eq:ExactNumeratorAsymtpt.Formula.9}
\begin{array}{l}
\displaystyle
A(\theta_s)\,\sim\,\left(\dfrac {\theta_s}2\right)^\eta\,\sim\,2^{-\eta/2} z^{\eta/4} k^{-\eta/2},\qquad\qquad k\to\infty,
\end{array}
\end{equation}
and, from Eq.~(\ref{App:Eq:ExactNumeratorAsymtpt.Formula.4}),
\begin{equation}
  \label{App:Eq:ExactNumeratorAsymtpt.Formula.10}
\begin{array}{l}
\displaystyle
\phi_+''(\theta_s)\,=\,\dfrac{3\,+\,\cos\theta_s}{8\,\sin^3\dfrac{\theta_s} 2}\,{\sqrt z} \sim
\dfrac{4{\sqrt z}}{\theta_s^3}\,\sim\ 2^{1/2}\,z^{-1/4}\,k^{3/2},\qquad\qquad k\to \infty.
\end{array}
\end{equation}
so that, on evaluating at $\theta_s$,
\begin{equation}
\label{App:Eq:ExactNumeratorAsymtpt.Formula.11}
\begin{array}{l}
\displaystyle
\phi_+''(\theta_s)\,\sim\,
\dfrac{4{\sqrt z}}{\theta_s^3}\,\sim\ 2^{1/2}\,z^{-1/4}\,k^{3/2},\qquad\qquad k\to \infty,
\end{array}
\end{equation}
and
\begin{equation}
\label{App:Eq:ExactNumeratorAsymtpt.Formula.12}
\begin{array}{l}
\displaystyle
\phi(\theta_s)+\dfrac\pi 4\,\mathrm{sign}\{\phi''(\theta_s)\}\,\sim\,
2\sqrt 2\,z^{1/4}\,k^{1/2}\,-\,\sqrt z\,+\,\dfrac\pi 4,
\qquad\qquad k\to \infty.
\end{array}
\end{equation}
Finally, on substituting from Eqs. (\ref{App:Eq:ExactNumeratorAsymtpt.Formula.9}) - (\ref{App:Eq:ExactNumeratorAsymtpt.Formula.12}) into Eq.~(\ref{App:Eq:ExactNumeratorAsymtpt.Formula.8}),
simple algebra eventually leads to 
\begin{equation}
  \label{App:Eq:ExactNumeratorAsymtpt.1}
\begin{array}{l}
\displaystyle
\mathcal{I}_+\sim
\dfrac{\sqrt\pi\, z^{\frac{2+q}8}}{(2k)^{1+\frac q4}}\,\sin\left(2\sqrt{2} z^{1/4} k^{1/2}-\sqrt{z}-\dfrac{q\pi}4\right)
,\qquad\qquad k\to\infty.
\end{array}
 \end{equation}

Consider now the second contribution to the integral, namely $\mathcal{I}_-$.
There are no stationary points within the domain $\theta \in [0,\pi]$, since the real function
\begin{equation}
  \label{App:Eq:ExactNumeratorAsymtpt:IMeno.1}
\begin{array}{l}
\displaystyle
\phi'_-(\theta)\,=\,
k+\dfrac{\sqrt z}2\,\dfrac {\cos\dfrac\theta 2}{\sin^2\dfrac\theta 2}\,, 
\end{array}
 \end{equation}
turns out to be strictly positive within the whole integration domain. 
It is worth recasting the integral $\mathcal{I}_-$ as follows:
\begin{equation}
  \label{App:Eq:ExactNumeratorAsymtpt:IMeno.2}
\begin{array}{l}
\displaystyle
\mathcal{I}_-\,=\,
\dfrac 12\,
\mathrm{Im}\left\{
\exp\left(-\mathrm{i}\,\eta\dfrac\pi 2\right)
\int_0^\pi\,\mathrm{d}\theta\,
A(\theta)\,
\exp\left[\mathrm{i}\phi_-(\theta)\right]
\right\},
\end{array}
\end{equation}
and trying to estimate the truncated Fourier transform on using standard procedures, like those outlined for example in~\cite[Ch. 6]{Bender/Orszag/1978}. In particular, due the absence of stationary points, it would be expected the most important contributions to the integral coming from the endpoints $\theta=0$ and $\theta=\pi$. To this end, partial integration will be used. 
First of all, we have~\cite[Sec. 6.3]{Bender/Orszag/1978}
\begin{equation}
  \label{App:Eq:ExactNumeratorAsymtpt:IMeno.3}
\begin{array}{l}
\displaystyle
\int_0^\pi\,\mathrm{d}\theta\,
A(\theta)\,
\exp\left[\mathrm{i}\phi_-(\theta)\right]\,=\,
\int_0^\pi\,\mathrm{d}\theta\,
\dfrac{A(\theta)}{\mathrm{i}\phi'_-(\theta)}\,
\dfrac{\mathrm{d}}{\mathrm{d}\theta}\,
\exp\left[\mathrm{i}\phi_-(\theta)\right]\,=\,\\
\\
\,=\,
\displaystyle
\left[
\dfrac{A(\theta)}{\mathrm{i}\phi'_-(\theta)}
\right]^\pi_0\,-\,
\int_0^\pi\,\mathrm{d}\theta\,
\exp\left[\mathrm{i}\phi_-(\theta)\right]
\dfrac{\mathrm{d}}{\mathrm{d}\theta}
\left(\dfrac{A(\theta)}{\mathrm{i}\phi'_-(\theta)}\right),
\end{array}
\end{equation}
where
\begin{equation}
  \label{App:Eq:ExactNumeratorAsymtpt:IMeno.6}
\left\{
\begin{array}{l}
\displaystyle
\dfrac{A(\theta)}{\phi'_-(\theta)}\,\sim\,\dfrac{1}{2^{\eta+1}\sqrt z}\,\theta^{\eta+2}\,,\qquad \theta\to 0^+,\\
\\
\dfrac{A(\theta)}{\phi'_-(\theta)}\,\sim\,\dfrac{1}{2\sqrt 2\,\,k}\,(\pi-\theta)^{3/2}\,,\qquad \theta\to \pi^-,
\end{array}
\right.
\end{equation}
so that 
\begin{equation}
  \label{App:Eq:ExactNumeratorAsymtpt:IMeno.7}
\begin{array}{l}
\displaystyle
\int_0^\pi\,\mathrm{d}\theta\,
A(\theta)\,
\exp\left[\mathrm{i}\phi_-(\theta)\right]\,=\,
\int_0^\pi\,\mathrm{d}\theta\,
A_1(\theta)\,\exp\left[\mathrm{i}\phi_-(\theta)\right],
\end{array}
\end{equation}
where
\begin{equation}
  \label{App:Eq:ExactNumeratorAsymtpt:IMeno.7}
\begin{array}{l}
\displaystyle
A_1(\theta)\,=\,
-\dfrac{\mathrm{d}}{\mathrm{d}\theta}
\left(\dfrac{A(\theta)}{\mathrm{i}\phi'_-(\theta)}\right).
\end{array}
\end{equation}
On repeating the same procedure for $A_1(\theta)$
we have
\begin{equation}
  \label{App:Eq:ExactNumeratorAsymtpt:IMeno.8}
\left\{
\begin{array}{l}
\displaystyle
\dfrac{A_1(\theta)}{\phi'_-(\theta)}\,\sim\,\dfrac{(2+\eta)\,\mathrm{i}}{2^{\eta+2} \sqrt z}\,\theta^{\eta+3}\,,\qquad \theta\to 0^+,\\
\\
\dfrac{A_1(\theta)}{\phi'_-(\theta)}\,\sim\,-\dfrac{3\mathrm{i}}{4\,\sqrt 2\,\,k^2}\,(\pi-\theta)^{1/2}\,,\qquad \theta\to \pi^-,
\end{array}
\right.
\end{equation}
and then,
\begin{equation}
  \label{App:Eq:ExactNumeratorAsymtpt:IMeno.9}
\begin{array}{l}
\displaystyle
\int_0^\pi\,\mathrm{d}\theta\,
A(\theta)\,
\exp\left[\mathrm{i}\phi_-(\theta)\right]\,=\,
\int_0^\pi\,\mathrm{d}\theta\,
A_2(\theta)\,\exp\left[\mathrm{i}\phi_-(\theta)\right],
\end{array}
\end{equation}
where
\begin{equation}
  \label{App:Eq:ExactNumeratorAsymtpt:IMeno.10}
\begin{array}{l}
\displaystyle
A_2(\theta)\,=\,
-\dfrac{\mathrm{d}}{\mathrm{d}\theta}
\left(\dfrac{A_1(\theta)}{\mathrm{i}\phi'_-(\theta)}\right)\,=\,
-\dfrac{\mathrm{d}}{\mathrm{d}\theta}
\left[
\dfrac{1}{\phi'_-(\theta)}
\dfrac{\mathrm{d}}{\mathrm{d}\theta}\,
\left(
\dfrac{A(\theta)}{\phi'_-(\theta)}
\right)
\right].
\end{array}
\end{equation}
It should be  noted that now
\begin{equation}
  \label{App:Eq:ExactNumeratorAsymtpt:IMeno.12}
\left\{
\begin{array}{l}
\displaystyle
{A_2(\theta)}\,\sim\,-\dfrac{(\eta^2+5\eta+6)}{2^{\eta+2}\,\, z}\,\theta^{\eta+2}\,,\qquad \theta\to 0^+,\\
\\
{A_2(\theta)}\,\sim\,-\dfrac{3\,k^{-2}}{8\,\sqrt 2}\,(\pi-\theta)^{-1/2}\,,\qquad \theta\to \pi^-.
\end{array}
\right.
\end{equation}
from which a divergence at the endpoint $\theta=\pi$ appears. Then, partial integration is stopped and the 
integral in the right side of  Eq.~(\ref{App:Eq:ExactNumeratorAsymtpt:IMeno.9}) can be estimated via Watson's lemma~\cite[Sec. 6.4]{Bender/Orszag/1978}, which gives
%
\begin{equation}
  \label{App:Eq:ExactNumeratorAsymtpt:IMeno.15}
\begin{array}{l}
\displaystyle
\int_{0}^\pi\,\mathrm{d}\theta\,
A_2(\theta)\,
\exp\left[\mathrm{i}\phi_-(\theta)\right] \sim
(-1)^{k+1}\,\exp\left(-\mathrm{i}\dfrac \pi 4\right)\,\dfrac{3\,\sqrt{2\pi}}{16}\,k^{-5/2}\,,\qquad\qquad k \to \infty.
\end{array}
\end{equation}
Finally, on substituting from  Eq.~(\ref{App:Eq:ExactNumeratorAsymtpt:IMeno.15}) into Eq.~(\ref{App:Eq:ExactNumeratorAsymtpt:IMeno.2}), we arrive at
%
\begin{equation}
  \label{App:Eq:ExactNumeratorAsymtpt:IMeno.16}
\begin{array}{l}
\displaystyle
\mathcal{I}_-\,\sim\,
(-1)^k\,\dfrac{3\sqrt{2\pi}}{32\,k^{5/2}}\,\sin\left(\dfrac{q+2}4\,\pi\right), \qquad\qquad\qquad\qquad\qquad  k\to\infty.
\end{array}
\end{equation}
The comparison between Eqs. (\ref{App:Eq:ExactNumeratorAsymtpt.1}) and (\ref{App:Eq:ExactNumeratorAsymtpt:IMeno.16})
shows that the leading-order asymptotic behavior of the integral is dominated by the saddle-point contribution $\mathcal{I}_+$.
Specifically, we have
\begin{equation}
\dfrac{\mathcal{I}_-}{\mathcal{I}_+}=\mathcal{O}\left(k^{-\frac{q+6}4}\right),
\end{equation}
and since $q\in (-1,1)$, the exponent  $-(6+q)/4$ remains strictly negative, thus ensuring that ${\mathcal{I}_-}$ decays faster
than ${\mathcal{I}_+}$ for all values of $q$. 
As a consequence, the former  can safely be considered asymptotically negligible with respect to $\mathcal{I}_+$, thus proving  Eq.~(\ref{Eq:ExactNumeratorAsymtpt.1}).
}

\section{Wynn's $\epsilon$ algorithm}
\label{App:WynnEpsAlg}
%

The modern era of sequence transformations started with two seminal
articles by~\citet{Shanks/1955} and~\citet{Wynn/1956a},
respectively. Shanks introduced a powerful sequence transformation that computes Pad\'e approximants,
while Wynn showed that Pad\'{e} approximants can be computed effectively by means of the so-called epsilon algorithm,
corresponding  to the following nonlinear recursive
scheme~\citep[Eq.\ (4)]{Wynn/1956a}:
\begin{subequations}
  \label{eps_al}
  \begin{align}
    \label{eps_al_a}
    \varepsilon_{-1}^{(n)} & \; = \; 0 \, , 
    \qquad \varepsilon_0^{(n)} \, = \,
    s_n \, ,
    \qquad  n \in \mathbb{N}_{0} \, , \\
    \label{eps_al_b}
    \varepsilon_{k+1}^{(n)} & 
    \; = \; \varepsilon_{k-1}^{(n+1)} \, + \,
    \frac{1}{\varepsilon_{k}^{(n+1)} - \varepsilon_{k}^{(n)}} \, , 
    \qquad k, n \in \mathbb{N}_{0} \, .
  \end{align}
\end{subequations}
The elements $\varepsilon_{2k}^{(n)}$ with \emph{even} subscripts provide
approximations to the (generalized) limit $s$ of the input sequence $\{
s_{n} \}_{n=0}^{\infty}$, whereas the elements $\varepsilon_{2k+1}^{(n)}$
with \emph{odd} subscripts are only auxiliary quantities which diverge if
the whole process converges. 

If the elements of the input sequence $\{ s_{n} \}_{n=0}^{\infty}$ were
the partial sums of a (formal) power series 
$f(z)$, then the $\epsilon$ algorithm produces Pad\'{e} approximants to 
$f(z)$ according to
\begin{equation}
  \label{Eps_Pade}
  \varepsilon_{2 k}^{(n)} \; = \; [ k + n / k ]_f (z) \, ,
  \qquad k, n \in \mathbb{N}_{0} \, .
\end{equation}
The epsilon algorithm is not restricted to input data that are the
partial sums of a (formal) power series.  Therefore, it is more general
and more widely applicable than Pad\'{e} approximants. As a recent review
we recommend~\citep*{Graves-Morris/Roberts/Salam/2000a}.


\end{document}